\documentclass[12pt,a4paper,left,final]{article}
\usepackage[utf8]{inputenc}
\usepackage{amsfonts}
\usepackage{amssymb}
\usepackage{amscd,amssymb,stmaryrd}
\usepackage[fleqn]{amsmath}
\usepackage{fancybox}
\usepackage{amscd,amstext}
\usepackage{hyperref}
\usepackage{mathabx}
\usepackage{color}
\usepackage{caption}
\usepackage{multicol}
\usepackage{cite}
\usepackage{amsthm}
\usepackage{mathtools}
\usepackage{bigints}
\usepackage{amsrefs}
\usepackage{showlabels} 
\usepackage{makeidx}
\usepackage{enumitem}
\usepackage{mathrsfs}
\usepackage{tikz}
\usepackage{empheq}
\usepackage{framed}
\usepackage{ulem}
\usepackage{nameref}

\numberwithin{equation}{section}

\newcommand{\ds}{\displaystyle}
\newcommand{\ra}{\rightarrow}

\newtheorem{theorem}{Theorem}
\newtheorem{lemma}{Lemma}

\newtheorem{definition}{Definition}

\newtheorem{assumption}{Assumption}
\newtheorem*{theorem*}{Theorem}
\newtheorem*{lemma*}{Lemma}
\newtheorem*{conj*}{Conjecture}
\newtheorem*{corollary*}{Corollary}
\newtheorem*{proposition*}{Proposition}
\newtheorem*{example*}{Example}
\newcommand{\rom}[1]{\uppercase\expandafter{\romannumeral #1\relax}}

\newcommand{\ul}{\uline}

\newcommand{\lb}{\lbrace}
\newcommand{\rb}{\rbrace}
\newcommand{\lan}{\langle}
\newcommand{\ran}{\rangle}

\newcommand{\R}{\mathbb{R}}

\newcommand{\N}{\mathbb{N}}
\newcommand{\E}{\mathbb{E}}

\newcommand{\prob}{\mathbb{P}}

\newcommand{\Holder}{H\"older}

\newcommand{\cS}{\mathcal{S}}

\newcommand{\cV}{\mathcal{V}}

\newcommand{\cO}{\mathcal{O}}

\newcommand{\sB}{\mathscr{B}}

\newcommand{\sF}{\mathscr{F}}

\newcommand{\sL}{\mathscr{L}}
\newcommand{\sA}{\mathscr{A}}

\newcommand{\Ga}{\alpha}
\newcommand{\Gb}{\beta}
\newcommand{\Gd}{\delta}
\newcommand{\Ge}{\varepsilon}
\newcommand{\Gg}{\gamma}

\newcommand{\Gl}{\lambda}
\newcommand{\GL}{\Lambda}

\newcommand{\Gs}{\sigma}
\newcommand{\Gt}{\tau}
\newcommand{\Gth}{\theta}
\newcommand{\GG}{\Gamma}
\newcommand{\GO}{\Omega}
\newcommand{\GD}{\Delta}

\newcommand{\Pot}{Proof of Theorem }

\theoremstyle{remark}

\makeatletter
\newcommand{\labitem}[2]{%
\def\@itemlabel{\textbf{#1}}
\item
\def\@currentlabel{#1}\label{#2}}
\makeatother
\date{}
\begin{document}

\title{\bf Periodic measures for a class of SPDEs with regime-switching}

\maketitle

\vskip -0.7cm

\centerline{Chun Ho Lau and Wei Sun}
\centerline{\small  Department of
Mathematics and Statistics}
\centerline{\small  Concordia University}
\centerline{\small Montreal, H3G 1M8, Canada}
\centerline{\small
E-mail: chunho.lau@concordia.ca, wei.sun@concordia.ca}

\vskip 1.7cm

\abstract{}

\noindent We use the variational approach to investigate periodic measures for a class of SPDEs with regime-switching. The hybrid system is driven by degenerate L\'{e}vy noise. We use the Lyapunov function method to study the existence of periodic measures and show the uniqueness of periodic measures by establishing the strong Feller property and irreducibility of the associated time-inhomogeneous semigroup. The main results are applied to stochastic porous media equations with regime-switching.

\endabstract

\vskip 0.5cm


\vskip 0.3cm


\vskip 1.5cm

\section{Introduction}


Stochastic partial differential equations (SPDEs) have gained more and more attentions in recent years. They have been applied to different fields including physics, biology, economics, etc. See  \cite{DZ3}, \cite{PZ} and \cite{LR} for general discussions on SPDEs and their applications.

\vskip 0.1cm

It is important to  study long time behaviours of solutions to SPDEs. We refer the reader to Da Prato and Zabczyk \cite{DZ} and Maslowski and Seidler \cite{MS} for systematic investigation of ergodicity for time-homogeneous SPDEs. In the past decades, many new results have been obtained for the existence and uniqueness of invariant measures.  Here we list some of them which motivated our paper. Hairer and Mattingly \cite{HM}  established ergodicity of the 2D Navier-Stokes equations with degenerate stochastic
forcing. Romito and Xu \cite{RX} discussed invariant measures of the 3D stochastic Navier-Stokes equations driven by mildly degenerate noise. Xie \cite{X} obtained the uniqueness of invariant measures for general SPDEs driven by non-degenerate L\'{e}vy noise. Wang \cite{W} used Harnack inequalities to investigate ergodicity of SPDEs. Gess and R\"{o}ckner \cite{GR} studied regularity and characterization for quasilinear SPDEs driven by degenerated Wiener noise. Zhang \cite{Zr} considered invariant measures of 3D stochastic MHD-$\Ga$ model driven by degenerate noise.  Neu\ss \ \cite{Neuss} studied ergodicity for singular-degenerate stochastic porous media equations.

\vskip 0.1cm

If stochastic equations are time-inhomogeneous, in general, we do not expect that invariant measures exist; instead, we consider periodic measures. There are many existing results studying the periodic behavior of stochastic differential equations (SDEs) and SPDEs.  In \cite{K}, Khasminskii systematically studied periodically varying properties of SDEs driven by Brownian motions. In \cite{ZWL}, Zhang et al.
investigated the existence and uniqueness of periodic solutions of SDEs driven by L\'evy processes.
In \cite{GS}, Guo and Sun generalized Doob's celebrated theorem on the uniqueness of invariant measures for time-homogeneous Markov processes so as to obtain the ergodicity and
uniqueness of periodic solutions for non-autonomous SDEs driven by L\'evy noises. For some other results related to periodic measures of SDEs, we refer the reader to Da Prato and Tudor \cite{DT}, Xu et al. \cite{a9,ra}, Chen et al. \cite{CHLY}, Hu and Xu \cite{a12}, and Ji et al. \cite{JQSY}.
In \cite{FZ}, Feng and Zhao
showed that there exist pathwise random periodic solutions to some SPDEs. In \cite{DD}, Da Prato and Debussche investigated the long time behavior of solutions to the 2D Stochastic Navier-Stokes equations with a time-periodic forcing term. In \cite{CL}, Cheng and Liu  used the variational approach to study recurrent properties
of solutions to SPDEs driven by Wiener noise. Under suitable conditions, in particular, by assuming strict monotonicity, they showed that the recurrent solutions are globally asymptotically stable in square-mean sense. In \cite{YB},  Yuan and Bao used the semigroup method to establish the exponential stability for a class of  finite regime-switching SPDEs driven by L\'evy noise.

\vskip 0.1cm

The aim of this work is to investigate the existence and uniqueness of periodic measures for a class of SPDEs with regime switching. The model consists of two
component processes $(X(t), \Lambda(t))$ with $X(t)$ and $\Lambda(t)$ being of continuous and
discrete states, respectively. The evolution of $X(t)$ is described by an SPDE that is driven by degenerate L\'evy noise. Through introducing regime switching $\Lambda(t)$ to the random dynamical system, more flexibility can be added in applications. The study of such a hybrid system is becoming more and
more important in different research areas such as biology, ecosystems, wireless
communications, signal processing, engineering and mathematical finance. Completely different from the methods of \cite{CL} and \cite{YB}, we will investigate ergodicity of  SPDEs with countable regime-switching through considering the strong Feller property and irreducibility of the corresponding time-inhomogeneous semigroups.

\vskip 0.1cm

Now we describe the framework of this paper. Let $(H, \lan\cdot,\cdot\ran_H)$ be a real separable Hilbert space and  $(V, |\cdot|_V)$ a real reflexive Banach space that is continuously and densely embedded into $H$. Denote by $V^*$ the dual space of $V$ and  $\lan \cdot, \cdot\ran$ the duality between $V$ and $V^*$. We have $ \lan u, v\ran = \lan u,v\ran_H$ for $u\in H$ and $v\in V$. Let $\{W(t)\}_{t\ge 0}$ be an $H$-valued cylindrical Wiener process on a complete filtered probability space $(\GO, \sF, \{\sF_t\}_{t\geq 0}, \prob)$. Denote by $L(H)$ and $L_2(H)$ the spaces of all bounded operators and Hilbert-Schmidt operators on $H$, respectively. Let $Z$ be a real Banach space with norm $|\cdot|$ and $N$ a Poisson random measure on $(Z,{\sB}(Z))$
with intensity measure $\nu$. We assume that $W$ and $N$ are independent. Set $\widetilde{N}(dt,dz)=N(dt,dz)-dt\nu(dz)$. Let $\cS=\lb 1,2,\dots\rb$.

\vskip 0.1cm

We consider the SPDE \begin{align} \label{maineqnR}
 dX(t)&= A(t,X(t),\GL(t))dt+B(t,X(t),\GL(t))dW(t) \nonumber\\
 &\quad\ +\int_{\{|z|<1\}} H(t,X(t),\GL(t),z)\widetilde{N}(dt,dz)\nonumber\\
 &\quad\  +\int_{\{|z|\ge1\}} J(t,X(t),\GL(t),z)N(dt,dz)
 \end{align}
with $X(0)=x\in H$. Hereafter, we assume that the coefficient functions $A:[0,\infty)\times V\times \cS \ra V^*$, $B:[0,\infty)\times V \times \cS\ra L_2(H)$ and $H,J: [0,\infty)\times V\times \cS \times Z \ra H$ are all measurable. The process $\GL(t)$ has state space $\cS$ such that when $\GD\ra 0$,
\begin{align*}
\prob(\GL(t+\GD)=j | \GL(t)=i, X(t)=x) =\begin{cases} q_{ij}(x)\GD+o(\GD),\quad &\text{if }i\neq j,\\ 1+ q_{ij}(x)\GD+o(\GD),  \quad &\text{if }i= j. \end{cases}
\end{align*}
Hereafter, $\{q_{ij}\}$ are Borel measurable functions on $H$ such that $q_{ij}(x)\geq 0$ for any $x\in H$ and $i,j\in\cS$ with $i\neq j$ and $\sum_{j\in\cS} q_{ij}(x)=0$ for any $x\in H$ and $i\in\cS$. Throughout this paper, we assume that
\begin{enumerate}
\labitem{(Q0)}{Q0}
$$L:=\sup_{x\in H, i\in\cS} \sum_{j\neq i} q_{ij}(x)<+\infty. $$
\end{enumerate}

We point out that $\GL(t)$ can be represented as a stochastic integral with
respect to a Poisson random measure. For each $x\in H$ and distinct $i,j\in\cS$, define $q_{i0}(x)=0$ and
$$\Delta_{ij}(x):= \bigg[\sum_{m=0}^{j-1} q_{im}(x), \sum_{m=0}^{j}q_{im} (x)\bigg).$$
Set
$$\GG(x,i,r)=\sum_{j\in \cS} (j-i)1_{\Delta_{ij}(x)}(r),\ \ \ \ (x,i,r)\in H\times {\cal S}\times [0,L].$$
Then, $\GL(t)$ can be modeled by
\begin{eqnarray}\label{222}
d\GL(t)=\int_{[0,L]}\GG(X(t-),\GL(t-),r) N_1(dt,dr)
\end{eqnarray}
for some  Poisson random measure $N_1$ with the Lebesgue measure on $[0,L]$ as its characteristic measure. We assume that
$N_1(\cdot,\cdot)$ is independent of $W(\cdot)$ and  $N(\cdot,
\cdot)$.

\vskip 0.1cm

The remainder of this paper is  as follows. First, we show the existence and uniqueness of solutions to the hybrid system (\ref{maineqnR}) and (\ref{222}) in Section 2. Then, in Sections 3, we establish the strong Feller property and irreducibility for the time-inhomogeneous semigroup corresponding to the hybrid system. We adopt the coupling method used in Zhang \cite{Z}. This remarkable  coupling method by change of measure was first introduced  to establish the dimension-free Harnack inequality by Arnaudon et al. \cite{ATW}. In Section 4, we obtain the existence and uniqueness of periodic measures for  the hybrid system by using the Lyapunov function method and generalizing the method of Guo and Sun \cite{GS} from SDEs to SPDEs. Finally, in Section 5, we use stochastic porous media equations as an example to illustrate the theory. The main theorem of this  paper is Theorem \ref{existencePSRS}.


\section{Existence and Uniqueness of Solutions}
To establish the existence and uniqueness of solutions to Equations \eqref{maineqnR} and \eqref{222}, we impose the following assumption.

\begin{assumption} \label{A1}
Suppose that there exist $\Ga>1$, $\Gb \geq 0$, $\Gth>0$, $K\in \R$, $\Gg<\frac{\Gth}{2\Gb}$, $c>0$, $\rho\in L^{\infty}_{loc}(V;[0,\infty))$ and $C\in L^{\frac{\Gb+2}{2}}_{loc}([0,\infty);[0,\infty))$ such that for $v_1,v_2,v\in V$, $i\in\cS$ and $t\in [0,\infty)$,

\begin{enumerate}[leftmargin=\parindent,align=left,labelwidth=\parindent]
\labitem{(HC)}{hemicontr} (Hemicontinuity) $s\mapsto \lan A(t,v_1+sv_2,i),v\ran$ is continuous on $\R$.

\labitem{(LM1)}{localmonr} (Local monotonicity)
\begin{eqnarray*}
&&2\lan A(t,v_1,i)-A(t,v_2,i),v_1-v_2\ran + \Vert B(t,v_1,i)-B(t,v_2,i)\Vert_{L_2(H)}^2\nonumber\\
&&+\int_{\{|z|<1\}} |H(t,v_1,i,z)-H(t,v_2,i,z)|_H^2 \nu(dz) \nonumber\\
&\leq& [K+\rho(v_2)]|v_1-v_2|_H^2.
\end{eqnarray*}

\labitem{(C)}{coercivityr} (Coercivity)
\begin{eqnarray*}
&&2\lan A(t,v,i),v\ran+ \Vert B(t,v,i)\Vert_{L_2(H)}^2 +\int_{\{|z|<1\}} |H(t,v,i,z)|_H^2 \nu(dz)\nonumber\\
& \leq& C(t)-\Gth |v|_V^{\Ga}+c|v|_H^2.
\end{eqnarray*}

\labitem{(G1)}{growthAr}(Growth of $A$)
\begin{align*}
|A(t,v,i)|_{V^*}^{\frac{\Ga}{\Ga-1}} \leq [C(t)+c|v|_V^{\Ga}](1+|v|_H^{\Gb}).
\end{align*}
\labitem{(G2)}{growthBHr}(Growth of $B$ and $H$)
\begin{align*}
\Vert B(t,v,i)\Vert_{L_2(H)}^2+\int_{\{|z|<1\}} |H(t,v,i,z)|_H^2 \nu(dz)  \leq C(t)+\Gg|v|_V^{\Ga}+c|v|_H^{2}.
\end{align*}
\labitem{(G$\Gb$)}{growthHr}(Growth of $H$ in $L^{\Gb+2}$)
\begin{align*}
\int_{\{|z|<1\}} |H(t,v,i,z)|_H^{\Gb+2} \nu(dz)  \leq [C(t)]^{\frac{\Gb+2}{2}}+c|v|_H^{\Gb+2}.
\end{align*}
\labitem{(G$\rho$)}{growthRho} (Growth of $\rho$)
\begin{align*}
\rho(v)\leq c(1+|v|_V^{\Ga})(1+|v|_H^{\Gb}).
\end{align*}
\end{enumerate}
\end{assumption}

Now we can state the main result of this section.
\begin{theorem} \label{2-1}
Suppose that Assumption \ref{A1} and condition \ref{Q0} hold. Let $T>0$, $x\in H$ and $i\in\cS$. Then, there exists a unique $H\times \cS$-valued adapted c\`adl\`ag process $\lb (X(t),\GL(t))\rb_{t\in[0,T]}$ such that
\begin{enumerate}
\item any $dt\times \prob$-equivalent class $\widehat{X}$ of $X$ is in $L^{\Ga}([0,T];V)\bigcap L^2([0,T];H)$, $\prob$-a.s.;
\item for any $V$-valued progressively measurable $dt\times \prob$-version $\overline{X}$ of $\widehat{X}$, the following holds for all $t\in[0,T]$ and $\prob$-a.s.:
\begin{align}\label{May14}
 X(t)&= x+\int_0^t A(s,\overline{X}(s),{\GL}(s))ds+\int_0^t B(s,\overline{X}(s),{\GL}(s))dW(s) \nonumber\\
 &\quad\ +\int_0^t \int_{\{|z|<1\}} H(s,\overline{X}(s),{\GL}(s),z)\widetilde{N}(ds,dz)\nonumber\\
&\quad\ +\int_0^t\int_{\{|z|\ge1\}} J(s,\overline{X}(s),{\GL}(s),z)N(ds,dz);
\end{align}
\item $\Lambda(0)=i$ and Equation  \eqref{222} holds.

\end{enumerate}
\end{theorem}

\vskip 0.2cm

\begin{proof}\ \ First, we consider Equation \eqref{maineqnR} when $i\in\cS$ is fixed. For this case,  \eqref{maineqnR} becomes an ordinary SPDE.
To simplify notation, we drop the dependence on $i$. We have the following result on the existence and uniqueness of solutions.

\vskip 0.4cm

\begin{theorem}{(\cite{BLZ}*{Theorem 1.2})} \label{2-2}\ \
Under the assumptions of Theorem \ref{2-1}, there exists a unique $H$-valued adapted c\`adl\`ag process $\lb X(t)\rb_{t\in[0,T]}$ such that
\begin{enumerate}
\item any $dt\times \prob$-equivalent class $\widehat{X}$ of $X$ is in $L^{\Ga}([0,T];V)\bigcap L^2([0,T];H)$, $\prob$-a.s.;
\item for any $V$-valued progressively measurable $dt\times \prob$-version $\overline{X}$ of $\widehat{X}$, the following holds for all $t\in[0,T]$ and $\prob$-a.s.:
\begin{align*}
 X(t)&= x+\int_0^t A(s,\overline{X}(s))ds+\int_0^t B(s,\overline{X}(s))dW(s) \\
 &\quad\ +\int_0^t \int_{\{|z|<1\}} H(s,\overline{X}(s),z)\widetilde{N}(ds,dz)\\
&\quad\  +\int_0^t\int_{\{|z|\ge1\}} J(s,\overline{X}(s),z)N(ds,dz).
\end{align*}
\end{enumerate}
\end{theorem}

\vskip 0.2cm

By Theorem \ref{2-2}, we know that for any $(x,i)\in H\times\cS$, there exists a unique $H$-valued adapted process $X^{(i)}(t)$ such that
\begin{align} \label{RSThm1Eq1}
 X^{(i)}(t)&= x+\int_0^t A(s,\overline{X^{(i)}}(s),i)ds+\int_0^t B(s,\overline{X^{(i)}}(s),i)dW(s) \nonumber \\
 &\quad\ +\int_0^t \int_{\{|z|<1\}} H(s,\overline{X^{(i)}}(s),i,z)\widetilde{N}(ds,dz)\nonumber\\
&\quad\  +\int_0^t\int_{\{|z|\ge1\}} J(s,\overline{X^{(i)}}(s),i,z)N(ds,dz),
\end{align}
where $\overline{X^{(i)}}$ is a $V$-valued progressively measurable $dt\times \prob$-version.
Let $0<\Gs_1<\Gs_2<\cdots<\Gs_n<\cdots$ be the set of all jump points of the stationary point process $p_1(t)$ corresponding to the Poisson random measure $N_1(dt,dr)$. We have $\lim_{n\ra\infty}\Gs_n=\infty$ almost surely by condition \ref{Q0}.

\vskip 0.1cm

Now we construct the processes $(X,\GL)$ and its progressively measurable version.
For $t\in [0,\sigma_1)$, define
\begin{equation}\label{May14b}
(X(t),\GL(t))=(X^{(i)}(t), i),\ \ \ \ \overline{X}(t)= \overline{X^{(i)}}(t).
\end{equation}
Set
$$\GL(\sigma_1)=i+\sum_{j\in\cS}(j-i)1_{\Delta_{ij}(X^{(i)}(\sigma_1-))}(p_1(\sigma_1)). $$
Then, (\ref{May14}) holds for $t\in [0,\Gs_1)$.

\vskip 0.1cm

Let
$$
\widetilde{W}(t)=W(t+\sigma_1)-W(\sigma_1), \ \
\widetilde{p}(t)=p(t+\sigma_1),\ \
\widetilde{p_1}(t)=p_1(t+\sigma_1).$$ Set
\begin{eqnarray*}
  \begin{array}{l}
X^{(\Lambda(\sigma_1))}(0)=X^{(i)}(\sigma_1),\\
     (\widetilde{X}(t),\widetilde{\Lambda}(t)) =(X^{(\Lambda(\sigma_1))}(t),\Lambda(\sigma_1)),\ \ t\in[0,\sigma_2-\sigma_1), \\
          \widetilde{\Lambda}(\sigma_2-\sigma_1)=\Lambda(\sigma_1)+\sum\limits_{j\in\mathcal{S}}
  (j-\Lambda(\sigma_1))\mathrm{1}_{\widetilde{A}(j)}(\widetilde{p}_1(\sigma_2-\sigma_1)),
  \end{array}
\end{eqnarray*}
where
$$\widetilde{A}(j)=\Delta_{\Lambda(\sigma_1),j}(X^{(\Lambda(\sigma_1))}((\sigma_2-\sigma_1)-)).$$
Then, for $t\in [\Gs_1, \Gs_2)$,  we define
$$(X(t),\Lambda(t))=(\widetilde{X}(t-\sigma_1),\widetilde{\Lambda}(t-\sigma_1)),\ \ \ \ \overline{X}(t)= \overline{\widetilde{X}}(t).$$
which together with  (\ref{May14b}) gives the unique
solution on the time interval $[0, \sigma_2)$. Continuing
this procedure inductively, we define $(X(t), \Lambda(t))$ on the
time interval $[0, \sigma_n)$ for each $n$. Therefore, $(X(t),
\Lambda(t))$ is the unique (c\`adl\`ag) solution to the hybrid system
(\ref{maineqnR}) and (\ref{222}) since $\lim_{n\ra\infty}\Gs_n=\infty$ almost surely.
\end{proof}

\section{Ergodicity}

Let  $\lb (X(t),\GL(t))\rb_{t\geq 0}$ be the unique solution to the hybrid system (\ref{maineqnR}) and (\ref{222}). By the standard argument, we know that  $\lb (X(t),\GL(t))\rb_{t\geq 0}$ is a Markov process (cf. \cite[Theorem 4.8]{GM} and \cite[Theorem 6.4.5]{A}). Denote by $\sB(H\times \cS)$ the Borel $\sigma$-algebra of $H\times \cS$, and $B_b(H\times \cS)$ and $C_b(H\times \cS)$ the spaces of all real-valued bounded Borel measurable and continuous functions  on $H\times \cS$, respectively. Let $P(s,(x,i);t,A)$ be the transition probability function of $\lb (X(t),\GL(t))\rb_{t\geq 0}$ given by $$P(s,(x,i);t,A):= \prob\big((X(t),\GL(t))\in A|(X(s),\GL(s))=(x,i)\big), $$ where $x\in H, \ i\in \cS, \  A\in \sB(H\times \cS)$ and $0\leq s< t<\infty$. Define the corresponding time-inhomogeneous transition semigroup by
$$P_{s,t}f(x,i):= \int_{H\times \cS} f(w) P(s,(x,i);t,dw),\ \ \ \ f\in B_b(H\times \cS). $$
In this section, we will establish the strong Feller property and irreducibility of $\{P_{s,t}\}$.

\subsection{Strong Feller Property}

Denote by $L_2^{+,s}(H)$ the set of all Hilbert-Schmidt operators on $H$ that are positive and self-adjoint. We impose the following assumption.
\vskip 0.4cm
\begin{assumption} \label{A2}  Suppose that $\Ga\geq 2$ and the following conditions hold:

\begin{enumerate}
\item
\begin{enumerate}[leftmargin=\parindent,align=left,labelwidth=\parindent]
\labitem{(LipB)}{LipschitzBr} For any $n\in\N$, there exists $C_n>0$ such that
$$\Vert B(t,v_1,i)-B(t,v_2,i)\Vert_{L_2(H)}\leq C_n |v_1-v_2|_H $$
for all $v_1,v_2\in V$ with $|v_1|_H,|v_2|_H\leq n$, $t\ge 0$ and  $i\in\cS$.
\end{enumerate}

\item  There exist $\Gl\in [2,\infty)\cap(\Ga-2,\infty)$, $\lb B_n\rb\subset L_2^{+,s}(H)$ and $n_0\in \N$ such that the following conditions hold:
\begin{enumerate}[leftmargin=\parindent,align=left,labelwidth=\parindent]
\labitem{(N)}{nondegenr} For any $n \in \N$, $t\geq 0$, $v\in V$ with $|v|_H\leq n$ and $i\in\cS$,  $$B(t,v,i)[B(t,v,i)]^* \geq B_n^2.$$

\labitem{(LM2)}{Localmonotone2r}For any $n\geq n_0$, there exist $\widetilde{K_n}\geq 0$ and $\Gd_{n}>0$ such that
\begin{eqnarray*}
&&2\lan A(t,v_1,i)-A(t,v_2,i), v_1-v_2\ran + \Vert B(t,v_1,i)-B(t,v_2,i)\Vert_{L_2(H)}^2\\
&&+\int_{\{|z|<1\}} |H(t,v_1,i,z)-H(t,v_2,i,z)|_H^2\nu(du) \nonumber\\
&\leq& -\Gd_n |{B_n}^{-1}(v_1-v_2)|_H^{\Gl} |v_1-v_2|_H^{\Ga-\Gl}+\widetilde{K_n}|v_1-v_2|_H^2
\end{eqnarray*}
for all $v_1,v_2\in V$, $t\geq 0$ and $i\in\cS$.
\end{enumerate}
\end{enumerate}
\end{assumption}

\vskip 0.2cm

Now we state the main theorem of this subsection.

\vskip 0.4cm

\begin{theorem} \label{3-1}
Suppose that Assumption \ref{A1} holds with $C\in L^{\infty}_{loc}([0,\infty);(0,\infty))$, Assumption \ref{A2} and condition \ref{Q0} hold. Then, the transition semigroup $\lb P_{s,t}\rb$ of $(X(t),\GL(t))$ is strong Feller.
\end{theorem}

\vskip 0.1cm

We will first prove Theorem \ref{3-1} for the case that $\Lambda(t)\equiv i $ for some $i\in\cS$ in \S 3.1.1, and then give the proof for the general case in \S 3.1.2.

\subsubsection{Strong Feller Property for SPDEs}
Let $i\in \cS$ be fixed. Then,  we can treat Equation \eqref{maineqnR} as an ordinary SPDE. We will  generalize the remarkable method of Zhang \cite{Z} to include jumps. Our goal is to establish the following result.

\begin{theorem}\label{3-2}
Under the assumptions of Theorem \ref{3-1}, the transition semigroup $\lb P_{s,t}\rb$ of $X(t)$ is strong Feller.
\end{theorem}

To prove Theorem \ref{3-2}, we first consider the case that $J\equiv 0$ and establish a lemma that is similar to \cite{Z}. To simplify notation, we drop the dependence on $i$. Fix $T>0$.

\begin{lemma} \label{3-L-1}\ \ Suppose that Assumption \ref{A1} holds with $\Ga\geq 2$, $J\equiv 0$ and the following conditions hold:

\vskip 0.2cm

(i) There exists $K_2>0$ such that
$$|[B(t,v_1)-B(t,v_2)]^*(v_1-v_2)|_H \leq K_2 (|v_1-v_2|_H^2\wedge |v_1-v_2|_H) $$ for all $t\in [0,T]$ and $v_1,v_2\in V$.

\vskip 0.2cm

(ii) There exist $\Gl\in [2,+\infty)\cap (\Ga-2, +\infty)$, $\overline{B}\in L_2^{+,s}(H)$, $\Gd>0$, $\widetilde{K}\geq 0$ such that
$$B(t,v)[B(t,v)]^*\geq \overline{B}^2,
$$
and
\begin{eqnarray}\label{coereqn2}
&&2\lan A(t,v_1)-A(t,v_2), v_1-v_2\ran + \Vert B(t,v_1)-B(t,v_2)\Vert_{L_2(H)}^2\nonumber\\
&&+\int_{\{|z|<1\}} |H(t,v_1,z)-H(t,v_2,z)|_H^2\nu(du) \nonumber\\
&\leq& -\Gd |\overline{B}^{-1}(v_1-v_2)|_H^{\Gl} |v_1-v_2|_H^{\Ga-\Gl}+\widetilde{K}|v_1-v_2|_H^2
\end{eqnarray}
 for all $t\in [0,T]$ and $v,v_1,v_2\in V$.
\vskip 0.2cm

\noindent Then, $P_{s,t}f$ is $\frac{\Gl+2-\alpha}{2\Gl}$-\Holder\ continuous for any
$f\in B_b(H)$.
\end{lemma}

\vskip 0.2cm

\begin{proof} \ \ We follow the elegant method of \cite[Lemma 3.1]{Z}.
Let $\Ge\in (0,1)$ satisfying $0\vee (\Ga-2) <\Gl(1-\Ge)<(2\Ga-2)\wedge \Ga$. Take $\Ga'\in (0,\Ge)$, whose value will be determined at the end of the proof. For $x,y\in H$, consider
\begin{align*}
dX(t)=&A(t,X(t))dt+B(t,X(t))dW(t)+\int_{\{|z|<1\}} H(t,X(t),z)\widetilde{N}(dt,dz), \\
dY(t)=& A(t,Y(t))dt+B(t,Y(t))dW(t)+\int_{\{|z|<1\}} H(t,Y(t),z)\widetilde{N}(dt,dz)\\
&+|x-y|^{\Ga'}_H \frac{X(t)-Y(t)}{|X(t)-Y(t)|_H^{\Ge}} dt,
\end{align*}
with $X(0)=x$ and $Y(0)=y$, respectively. Define
$$\Gt_n:=\inf\left\{ t>0: |X(t)-Y(t)|_H\leq \frac{1}{n} \right\},\ \ n\in \mathbb{N},$$
and
$$\Gt := \lim_{n\rightarrow\infty}\tau_n.$$
By It\^{o}'s formula (cf. \cite{M}), we find that for $t<\tau$,
\begin{align}\label{21A}
d&|X(t)-Y(t)|_H^2\nonumber\\
&= 2\lan A(t,X(t))-A(t,Y(t)), X(t)-Y(t)\ran dt+\Vert B(t,X(t))-B(t,Y(t))\Vert_{L_2(H)}^2 dt \nonumber\\
&\quad + 2\lan X(t)-Y(t), [B(t,X(t))-B(t,Y(t))] dW(t)\ran_H \nonumber\\
&\quad +\int_{\{|z|<1\}}\bigg[|X(t)-Y(t)+H(t,X(t),z)-H(t,Y(t),z)|_H^2-|X(t)-Y(t)|_H^2 \bigg] \widetilde{N}(dt,dz)\nonumber \\
&\quad + \int_{\{|z|<1\}} |H(t,X(t),z)-H(t,Y(t),z)|_H^2 \nu(dz)dt \nonumber\\
&\quad -2|x-y|_H^{\Ga'}|X(t)-Y(t)|_H^{2-\Ge}dt,
\end{align}

\noindent which implies that
\begin{align*}
d&|X(t)-Y(t)|_H^2\\
&\leq  \left[-\Gd |\overline{B}^{-1}(X(t)-Y(t))|_H^{\Gl}|X(t)-Y(t)|_H^{\Ga-\Gl}\right.\\
&\left.\quad+\widetilde{K}|X(t)-Y(t)|_H^{2}-2|x-y|_H^{\Ga'}|X(t)-Y(t)|_H^{2-\Ge}\right]dt \\
&\quad + 2\lan X(t)-Y(t), [B(t,X(t))-B(t,Y(t))] dW(t)\ran_H \\
&\quad +\int_{\{|z|<1\}}\bigg[|X(t)-Y(t)+H(t,X(t),z)-H(t,Y(t),z)|_H^2-|X(t)-Y(t)|_H^2 \bigg] \widetilde{N}(dt,dz).
\end{align*}

\vskip 0.1cm
Set
$$r:= \frac{\Ga-\Gl(1-\Ge)}{2}.
$$
We have $r\in (0,1)$.
Set
$$V(x)=|x|^{1-r},\ \ x\neq 0.$$
We choose a sequence of functions $V_n\in C^2(\R)$ such that $V_n(x)=V(x)$ for $x\geq \frac{1}{n}$ and $V_n''(x)\leq 0$ for  $x\in\R$. Note that for $t<\tau_n$,
\begin{align*}
V_n'(|X(t)-Y(t)|_H^2)&=(1-r)|X(t)-Y(t)|_H^{-2r}, \\
V_n''(|X(t)-Y(t)|_H^2)&=-r(1-r)|X(t)-Y(t)|_H^{-2r-2}.
\end{align*}
Then, by (\ref{21A}), we obtain that for $t<\tau_n$,
{\small\begin{align}\label{21B}
d&|X(t)-Y(t)|_H^{2-2r}\nonumber \\
&=\bigg[ 2\lan A(t,X(t))-A(t,Y(t)), X(t)-Y(t)\ran +\Vert B(t,X(t))-B(t,Y(t))\Vert_{L_2(H)}^2\nonumber  \\ & \quad  +  \int_{\{|z|<1\}} |H(t,X(t),z)-H(t,Y(t),z)|_H^2 \nu(dz) -2|x-y|_H^{\Ga'}|X(t)-Y(t)|_H^{2-\Ge}\bigg]\nonumber\\
&\quad\cdot V_n'(|X(t)-Y(t)|_H^2)dt\nonumber  \\
&\quad+2V_n''(|X(t)-Y(t)|_H^2)\big|[B(t,X(t))-B(t,Y(t))]^*(X(t)-Y(t))\big|_H^2dt \nonumber\\ &\quad+ 2 V_n'(|X(t)-Y(t)|_H^2)  \lan X(t)-Y(t), [B(t,X(t))-B(t,Y(t))] dW(t)\ran_H \nonumber\\
&\quad +\int_{\{|z|<1\}} \Bigg[ V_n\bigg(|X(t)-Y(t)|_H^2\nonumber\\
&\quad\quad\quad + \Big[|X(t)-Y(t)+H(t,X(t),z)-H(t,Y(t),z)|_H^2-|X(t)-Y(t)|_H^2 \Big]\bigg)\nonumber\\ &\quad\quad\quad -V_n( |X(t)-Y(t)|_H^2)\Bigg]\widetilde{N}(dt,dz)\nonumber\\
&\quad +\int_{\{|z|<1\}} \Bigg\{ V_n\bigg(|X(t)-Y(t)|_H^2\nonumber\\
&\quad\quad\quad + \Big[|X(t)-Y(t)+H(t,z,X(t))-H(t,z,Y(t))|_H^2-|X(t)-Y(t)|_H^2 \Big]\bigg)\nonumber\\ &\quad\quad\quad -V_n'(|X(t)-Y(t)|_H^2) \Big[|X(t)-Y(t)+H(t,X(t),z)-H(t,Y(t),z)|_H^2-|X(t)-Y(t)|_H^2\Big]\nonumber\\ &\quad\quad\quad -V_n( |X(t)-Y(t)|_H^2) \Bigg\}\nu(dz)dt.
\end{align}}

By \eqref{coereqn2}, (\ref{21B}), the assumption on $B$ and the concave property of $V_n$, we obtain that for $t<\tau_n$,
{\small\begin{align} \label{2-2rIF2}
d&|X(t)-Y(t)|_H^{2-2r}\nonumber \\
&\leq -\Gd(1-r)|\overline{B}^{-1}(X(t)-Y(t))|_H^{\Gl}|X(t)-Y(t)|^{-\Gl\Ge}dt+\widetilde{K}(1-r)|X(t)-Y(t)|_H^{2-2r}dt\nonumber\\
&\quad  -2(1-r)|x-y|_H^{\Ga'}|X(t)-Y(t)|_H^{2-\Ge-2r}dt\nonumber\\
&\quad-2r(1-r)K^2_2\big(|X(t)-Y(t)|_H^{2-2r}\wedge |X(t)-Y(t)|_H^{-2r}\big)dt \nonumber\\
&\quad+ 2(1-r)|X(t)-Y(t)|_H^{-2r}\lan X(t)-Y(t), [B(t,X(t))-B(t,Y(t))] dW(t)\ran_H \nonumber\\
&\quad +\int_{\{|z|<1\}} \bigg[V_n(|X(t)-Y(t)+H(t,X(t),z)-H(t,Y(t),z)|_H^2)- V_n(|X(t)-Y(t)|_H^2)\bigg]\widetilde{N}(dt,dz).
\end{align}}

\noindent Then, by Gronwall's lemma, we get
\begin{align*}
\E\bigg[ \int_0^{T\wedge \tau_n} \frac{|\overline{B}^{-1}(X(t)-Y(t))|_H^{\Gl}}{|X(t)-Y(t)|_H^{\Gl\Ge}}dt \bigg]\leq \frac{ e^{(1-r)\widetilde{K}T}}{\Gd(1-r)} |x-y|_H^{2-2r}.
\end{align*}
Further, by Fatou's lemma, we get
\begin{align}\label{etannotbig2}
\E\bigg[ \int_0^{T\wedge \tau} \frac{|\overline{B}^{-1}(X(t)-Y(t))|_H^{\Gl}}{|X(t)-Y(t)|_H^{\Gl\Ge}}dt \bigg]\leq \frac{ e^{(1-r)\widetilde{K}T}}{\Gd(1-r)} |x-y|_H^{2-2r}.
\end{align}

\vskip 0.2cm

Define
$$\eta_n:= \inf \bigg\lb t>0: \int_0^{t} \frac{|\overline{B}^{-1}(X(s)-Y(s))|_H^{\Gl}}{|X(s)-Y(s)|_H^{\Gl\Ge}}ds \geq n \bigg\rb.$$
By \eqref{etannotbig2}, we get $\ds \lim_{n\ra\infty} \eta_n \geq T\wedge \tau. $ Set
$$\xi_n:= \eta_n\wedge \tau_n\wedge T.$$
Then,   we have that $\ds \lim_{n\ra\infty} \xi_n = T\wedge\tau$.
Define
{\small$$\widetilde{W}(t):= W(t)+\int_0^{t\wedge\tau} |x-y|_H^{\Ga'}\frac{[B(s,Y(s))]^*(B(s,Y(s))[B(s,Y(s))]^*)^{-1}(X(s)-Y(s))}{|X(s)-Y(s)|_H^{\Ge}}ds. $$}

\noindent Thus, $\lb \widetilde{W}(s)\rb_{s\in[0,\xi_n\wedge t]}$ is a cylindrical Wiener process on $H$ under the probability measure $R_{t\wedge\xi_n}\prob$. Define
$$E^{-1}(t,v):=[B(t,v)]^*(B(t,v))[B(t,v)]^*)^{-1},$$
and
\begin{align*}
R_t:&= \exp \bigg\lb -|x-y|_H^{\Ga'}\int_0^{t\wedge\tau} \bigg\lan \frac{E^{-1}(s,Y(s))(X(s)-Y(s))}{|X(s)-Y(s)|_H^{\Ge}}, dW(s)\bigg \ran_H \\
&\quad \quad \quad \quad \quad - \frac{|x-y|_{H}^{2\Ga'}}{2}\int_0^{t\wedge\tau} \frac{|E^{-1}(s,Y(s))(X(s)-Y(s))|_H^2}{|X(s)-Y(s)|_H^{2\Ge}} ds \bigg\rb.
\end{align*}
We will show that  $\lb \widetilde{W}(s)\rb_{s\in[0,T]}$ is a cylindrical Wiener process on $H$ under $R_{T}\prob$.

By \eqref{2-2rIF2} and the definition of $\widetilde{W}$, we obtain that for all $t<\tau_n$,
{\small\begin{align}\label{OP}
d&|X(t)-Y(t)|_H^{2-2r} \nonumber\\
&\leq  (1-r)\left[-\Gd\cdot\frac{|\overline{B}^{-1}(X(t)-Y(t))|_H^{\Gl}}{|X(t)-Y(t)|_H^{\Gl\Ge}}+\widetilde{K}|X(t)-Y(t)|_H^{2-2r}-2|x-y|_H^{\Ga'}|X(t)-Y(t)|_H^{2-2r-\Ge}\right]dt\nonumber\\
&\quad + 2(1-r)\left\lan \frac{[B(t,X(t))-B(t,Y(t))]^*(X(t)-Y(t))}{|X(t)-Y(t)|_H^{2r}},  d\widetilde{W}(t)\right\ran_H \nonumber\\
&\quad+ 2(1-r)|x-y|_H^{\Ga'}\bigg\lan \frac{X(t)-Y(t)}{|X(t)-Y(t)|_H^{2r}}, \frac{[B(t,X(t))-B(t,Y(t))]E^{-1}(t,Y(t))(X(t)-Y(t))}{|X(t)-Y(t)|_H^{\Ge}} \bigg\ran_H dt \nonumber\\
&\quad +\int_{\{|z|<1\}}\bigg[V_n(|X(t)-Y(t)+H(t,X(t),z)-H(t,Y(t),z)|_H^{2})-V_n(|X(t)-Y(t)|_H^{2}) \bigg] \widetilde{N}(dt,dz).
\end{align}}
By Young's inequality, we obtain that for $t<\tau$,
{\small\begin{align*}
& \bigg|2(1-r)|x-y|_H^{\Ga'}\bigg\lan \frac{X(t)-Y(t)}{|X(t)-Y(t)|_H^{2r}}, \frac{[B(t,X(t))-B(t,Y(t))]E^{-1}(s,Y(s))(X(s)-Y(s))}{|X(s)-Y(s)|_H^{\Ge}} \bigg\ran_H\bigg| \\
&\leq  \frac{\Gl-1}{\Gl}\cdot 2^{\frac{\Gl+1}{\Gl-1}}\Gd^{-\frac{1}{\Gl-1}}(1-r)|x-y|_H^{\frac{\Ga'\Gl}{\Gl-1}}K_2^{\frac{\Gl}{\Gl-1}}
\left(|X(t)-Y(t)|_H^{\frac{\Gl(2-2r)}{\Gl-1}}\wedge|X(t)-Y(t)|_H^{\frac{\Gl(1-2r)}{\Gl-1}}\right)\\
&\quad  + \frac{(1-r)\Gd}{2{\Gl}}\cdot\frac{|\overline{B}^{-1}[X(t)-Y(t)]|_{H}^{\Gl}}{|X(t)-Y(t)|_H^{\Gl\Ge}}.
\end{align*}}

\noindent Since $\frac{\Gl(1-2r)}{\Gl-1}\leq 2-2r$ and $|X(t)-Y(t)|_H^{\frac{\Gl(2-2r)}{\Gl-1}}\leq |X(t)-Y(t)|_H^{2-2r}$ if $|X(t)-Y(t)|_H\leq 1$, we get
{\small\begin{align*}
& \bigg|2(1-r)|x-y|_H^{\Ga'}\bigg\lan \frac{X(t)-Y(t)}{|X(t)-Y(t)|_H^{2r}}, \frac{[B(t,X(t))-B(t,Y(t))]E^{-1}(s,Y(s))(X(s)-Y(s))}{|X(s)-Y(s)|_H^{\Ge}} \bigg\ran_H\bigg| \\
&\leq   C_{\Gl,\Gd,K_2}(1-r)|x-y|_H^{\frac{\Ga'\Gl}{\Gl-1}}|X(t)-Y(t)|^{2-2r}_H + \frac{(1-r)\Gd}{2}\cdot\frac{|\overline{B}^{-1}[X(t)-Y(t)]|_{H}^{\Gl}}{|X(t)-Y(t)|_H^{\Gl\Ge}},
\end{align*}}

\noindent where  $C_{\Gl,\Gd,K_2}$ depends only on $\Gd$, $\Gl$ and $K_2$. Thus, for $t<\tau_n$, we have that
{\small\begin{align*}
d&|X(t)-Y(t)|_H^{2-2r} \\
&\leq  \bigg[-\frac{(1-r)\Gd}{2}\frac{|\overline{B}^{-1}(X(t)-Y(t))|_H^{\Gl}}{|X(t)-Y(t)|_H^{\Gl\Ge}}+(1-r)\widetilde{K}|X(t)-Y(t)|_H^{2-2r}\\
&\quad +C_{\Gl,\Gd,K_2}(1-r)|x-y|_H^{\frac{\Ga'\Gl}{\Gl-1}}|X(t)-Y(t)|^{2-2r}\bigg]dt\\
&\quad + 2(1-r)\left\lan \frac{[B(t,X(t))-B(t,Y(t))]^*(X(t)-Y(t))}{|X(t)-Y(t)|_H^{2r}},  d\widetilde{W}(t)\right\ran_H \\
&\quad +\int_{\{|z|<1\}}\bigg[V_n(|X(t)-Y(t)+H(t,X(t),z)-H(t,Y(t),z)|_H^{2})-V_n(|X(t)-Y(t)|_H^{2}) \bigg] \widetilde{N}(dt,dz).
\end{align*}

\noindent By  Gronwall's lemma, we get
{\small\begin{eqnarray}\label{EQSmain1}
&&\E_{R_{s\wedge\xi_n}\prob}\bigg(\int_{0}^{s\wedge\eta_n}\frac{|\overline{B}^{-1}(X(t)-Y(t))|_H^{\Gl}}{|X(t)-Y(t)|_H^{\Gl\Ge}}dt\bigg)\nonumber\\
&\leq& \frac{2\exp\left\{\left[C_{\Gl,\Gd,K_2}|x-y|_H^{\frac{\Ga'\Gl}{\Gl-1}}+\widetilde{K}\right](1-r)s\right\}}{(1-r)\Gd}|x-y|_H^{2-2r}.
\end{eqnarray}}}

By the definition of $\widetilde{W}$ and \eqref{EQSmain1}, we get
{\small\begin{eqnarray}
&&\sup_{\substack{s\in[0,T], n\in\N}}\E\big[R_{s\wedge\xi_n}\log(R_{s\wedge\xi_n})\big]\nonumber\\
&\leq& \sup_{\substack{s\in[0,T], n\in\N}} \frac{|x-y|_{H}^{2\Ga'}}{2}\E_{R_{s\wedge\xi_n}\prob}\bigg[\bigg(\int_{0}^{s\wedge \xi_n}\frac{|E^{-1}(t,Y(t))(X(t)-Y(t))|_H^{\Gl}}{|X(t)-Y(t)|_H^{\Gl\Ge}} dt\bigg)^{\frac{2}{\Gl}} T^{\frac{\Gl-2}{\Gl}} \bigg] \nonumber\\
&\leq& \left(\frac{T}{2}\right)^{\frac{\Gl-2}{\Gl}}\cdot \frac{\exp\left\{\left[C_{\Gl,\Gd,K_2}|x-y|_H^{\frac{\Ga'\Gl}{\Gl-1}}+\widetilde{K}\right]
(1-r)\frac{2T}{\Gl}\right\}}{[(1-r)\Gd]^{\frac{2}{\Gl}}}\cdot|x-y|_H^{\frac{4(1-r)}{\Gl}+2\Ga'}\label{RR}\\
& <&\infty.\nonumber
\end{eqnarray}}

\noindent Then $\lb R_{s\wedge \xi_n}\rb_{\substack{s\in[0,T], n\in\N}}$ and thus $\lb R_{s}\rb_{s\in[0,T]}$ is uniform integrable.  Hence, $\lb \widetilde{W}(t)\rb_{t\in[0,T]}$ is a cylindrical Wiener process under $R_T\prob$. Moreover, $\{Y(t)\}_{t\ge 0}$ satisfies
$$dY(t) = A(t,Y(t))dt+B(t,Y(t))d\widetilde{W}(t)+\int_{\{|z|<1\}} H(t,Y(t),z)\widetilde{N}(dt,dz)$$ with $Y(0)=y$. This implies that $\{Y(t)\}_{t\ge 0}$ is also a solution to \eqref{maineqnR} with $J\equiv 0$.

\vskip 0.1cm
We now show that $P_{s,t}f$ is $\frac{\Gl+2-\alpha}{2\Gl}$-\Holder\ continuous for any
$f\in B_b(H)$. To simplify notation, we only give the proof for the case that $s = 0$. The proof for the
case that $s > 0$ is completely similar.

\vskip 0.1cm

Let $f\in B_b(H)$, $0\leq t\leq T$ and $x,y\in H$. We have
\begin{align}\label{pp1}
|P_{0,t}f(x)-P_{0,t}f(y)| &= |\E[f(X(t))-R_tf(Y(t))]| \nonumber\\
&\leq\big|\E[f(Y(t))-R_tf(Y(t))]\big|+\big|\E\big\{[f(X(t))-f(Y(t))]1_{\lb\tau\geq t\rb}\big\}\big|\nonumber \\
&\leq |f|_{L^{\infty}} \big\{\E[|1-R_t|]+2\prob(\tau \geq t)\big\}.
\end{align}
By (\ref{etannotbig2}), (\ref{RR}) and the inequality
$$|1-e^x|\leq  xe^x+2|x|,\ \ x\in\R,$$
there exists $C'_{r,\widetilde{K},T,\Gd,\Gl,K_2}>0$ such that if $|x-y|_H\leq 1$,
\begin{align}\label{pp2}
\E[|1-R_t|]&\le \E[R_t\log R_t]+2\E[|\log R_t|]\nonumber\\
&\le  C'_{r,\widetilde{K},T,\Gd,\Gl,K_2} \bigg[|x-y|_H^{\frac{2(1-r)}{\Gl}+\Ga'}+|x-y|_H^{\frac{4(1-r)}{\Gl}+2\Ga'}\bigg].
\end{align}

Next, we estimate $\prob(\tau\geq t)$. Set $\widetilde{V}(x)=x^{\Ge/2}$ for $x\neq0$ and
choose a sequence of functions $\widetilde{V}_n\in C^2(\R)$ such that $\widetilde{V}_n(x)=\widetilde{V}(x)$ for $x\geq \frac{1}{n}$ and $\widetilde{V}''_n(x)\leq 0$ for  $x\in\R$. Similar to (\ref{OP}), we can show that for $s<\tau_n$,
{\small\begin{align*}
d&|X(s)-Y(s)|_H^{\Ge}\nonumber \\
&\leq \frac{\widetilde{K}\Ge}{2}|X(s)-Y(s)|_H^{\Ge}ds -\Ge|x-y|_H^{\Ga'}ds\nonumber  \\
&\quad+ 2\widetilde{V}_n'(|X(s)-Y(s)|_H^2)\lan X(s)-Y(s), [B(s,X(s))-B(s,Y(s))] dW(s)\ran_H \nonumber\\
&\quad +\int_{\{|z|<1\}} \bigg[ \widetilde{V}_n(|X(s)-Y(s)+H(s,X(s),z)-H(s,Y(s),z)|_H^2) -\widetilde{V}_n( |X(s)-Y(s)|_H^2)\bigg]\widetilde{N}(ds,dz),
\end{align*}}

\noindent which implies that
\begin{align} \label{expectation1}
\E\big[|X(s\wedge \tau_n)-Y(s\wedge \tau_n)|_H^{\Ge}\big]\leq& |x-y|_H^{\Ge}+\int_0^s\frac{\widetilde{K}\Ge}{2}\E\big[|X(u)-Y(u)|_H^{\Ge}\big]du\nonumber\\
&-\Ge |x-y|_{H}^{\Ga'} \E[s\wedge\tau_n].
\end{align}
By Gronwall's lemma, we get
$$\E\big[|X(s\wedge \tau_n)-Y(s\wedge \tau_n)|_H^{\Ge}\big]\leq e^{s\widetilde{K}\Ge/2} |x-y|_H^{\Ge}. $$
Then,
$$\int_0^t\E\big[|X(s\wedge \tau_n)-Y(s\wedge \tau_n)|_H^{\Ge}\big]ds\leq \frac{2}{\widetilde{K}\Ge} (e^{t\widetilde{K}\Ge/2}-1) |x-y|_H^{\Ge},$$
which together with \eqref{expectation1} implies that
$$\E[t\wedge\tau_n] \leq  \Ge^{-1}(e^{t\widetilde{K}\Ge/2}) |x-y|_H^{\Ge-\Ga'}.$$
Thus, we have that
\begin{equation}\label{pp3}
\prob(\tau> t)\leq \liminf_{n\ra\infty} \prob(\tau_n\geq t)\leq \liminf_{n\ra\infty} \frac{\E[t\wedge \tau_n]}{t} \leq \frac{e^{t\widetilde{K}\Ge/2}}{t\Ge}|x-y|_H^{\Ge-\Ga'}.
\end{equation}
Set
$$\Ga':=\frac{\Gl+2-\alpha}{2\Gl}.$$
Therefore, by (\ref{pp1}), (\ref{pp2}) and (\ref{pp3}), we obtain that for any $x,y\in H$ with $|x-y|_H<1$,
\begin{eqnarray*}
&&|P_{0,t}f(x)-P_{0,t}f(y)|\\
&\leq& |f|_{L^{\infty}}\left\{ C'_{r,\widetilde{K},T,\Gd,\Gl,K_2} \bigg[|x-y|_H^{\frac{2(1-r)}{\Gl}+\widetilde{\alpha}}+|x-y|_H^{\frac{4(1-r)}{\Gl}+2\widetilde{\alpha}}\bigg]+\frac{2e^{t\widetilde{K}\Ge/2}}{t\Ge}|x-y|_H^{\Ge-\widetilde{\alpha}}\right\}\\
&\leq& 2\left[C'_{r,\widetilde{K},T,\Gd,\Gl,K_2,\Ge}+\frac{e^{t\widetilde{K}\Ge/2}}{t\Ge}\right] |f|_{L^{\infty}} |x-y|^{\frac{\Gl+2-\alpha}{2\Gl}}_H.
\end{eqnarray*}
\end{proof}

\begin{proof}[\ul{\Pot \ref{3-2}}]\ \\

\noindent We first consider the case that $J\equiv 0$. For $R>0$, define
$$B_R(t,v):= \begin{cases} B(t,v) &\quad \text{if }|v|_H\leq R, \\ B(t,\frac{Rv}{|v|_H})&\quad \text{if }|v|_H> R.\end{cases} $$
Denote by $X(s,w;t)$ the solution of \eqref{maineqnR} with $X(s)=w$, $w\in H$ for fixed $i\in {\cal S}$. Suppose $|w|_H<R$. Define
\begin{align*}
\tau_R^w&:=\inf\lb t>s: |X(s,w;t)|_H\geq R \rb.
\end{align*}
Let $\{X_R(s,w;t)\}$ be the unique solution of the SPDE:
$$ dX_R(t)= A(t,X_R(t))dt+B_R(t,X_R(t))dW(t)+\int_{\{|z|<1\}} H(t,X_R(t),z)\widetilde{N}(dt,dz) $$ with $X(s)=w$. Denote by $\{P_{s,t}^R\}$ the  transition semigroup of $\{X_R(s,w;t)\}$. Suppose $x,y\in H$ with $|x|_H,|y|_H<R$. By the uniqueness of solutions, we find that $X(s,x;t)=X_R(s,x;t)$ and $X(s,y;t)=X_R(s,y;t)$ for all $t<\tau_R^x\wedge \tau_R^y$. Let $R>n_0$, which is given in condition \ref{Localmonotone2r}. By conditions \ref{nondegenr}, \ref{LipschitzBr}, \ref{Localmonotone2r} and replacing $\Gd$, $\widetilde{K}$, $\overline{B}$ with $\Gd_{{n_0}}$, ${K}_{\widetilde{n_0}}$,  $B_{{n_0}}$, respectively, we can apply Lemma \ref{3-L-1} to show that $\{P_{s,t}^R\}$ is strong Feller.

\vskip 0.1cm

Let $f\in B_b(H)$. We have
\begin{align}\label{AS1}
&|P_{s,t}f(x)-P_{s,t}f(y)|\nonumber\\
&\leq |\E[\{f(X(s,x;t))-f(X(s,y;t))\}1_{\lb\tau_R^x\wedge \tau_R^y> t\rb}]| + 2|f|_{L^{\infty}} \big[\prob(\tau_R^x\leq t)+\prob(\tau_R^y\leq t)\big]\nonumber \\
&= |\E[\{f(X_R(s,x;t))-f(X_R(s,y;t))\}1_{\lb\tau_R^x\wedge \tau_R^y> t\rb}]| + 2|f|_{L^{\infty}} \big[\prob(\tau_R^x\leq t)+\prob(\tau_R^y\leq t)\big] \nonumber\\
&\leq  |P_{s,t}^Rf(x)-P_{s,t}^Rf(y)|+ 2|f|_{L^{\infty}} \big[\prob(\tau_R^x\leq t)+\prob(\tau_R^y\leq t)\big].
\end{align}
By It\^{o}'s formula, we get
\begin{align*}
&d|X(s,x;t)|^2_H \\
&= \bigg[2\lan A(t,X(s,x;t)),X(s,x;t)\ran + \Vert B(t,X(s,x;t))\Vert_{L_2(H)}^2\\
&\quad +\int_{\{|z|<1\}} |H(t,X(s,x;t),z)|_H^2\nu(dz)\bigg]dt+ 2\lan B(t,X(s,x;t)),dW(t)\ran_H\\
&\quad + \int_{\{|z|<1\}} \big[|X(s,x;t)+H(t,X(s,x;t),z)|_H^2-|X(s,x;t)|_H^2\big]\widetilde{N}(dt,dz) \\
&\leq \sup_{u\in[0,t]}|C(u)|-\Gth|X(s,x;t)|_V^{\Ga}+c|X(s,x;t)|_H^2+ 2\lan B(t,X(s,x;t)),dW(t)\ran_H\\
&\quad + \int_{\{|z|<1\}} \big[|X(s,x;t)+H(t,X(s,x;t),z)|_H^2-|X(s,x;t)|_H^2\big]\widetilde{N}(dt,dz).
\end{align*}
Then, by Gronwall's lemma, we obtain that
\begin{align*}
\E\big[|X(s,x;t)|_H^2\big] \leq \left(\sup_{u\in[0,t]}|C(u)|+|x|_H^2\right)ce^{c(t-s)},
\end{align*}
and hence
\begin{align*}
\int_{s}^{t}\E[|X(s,x;u)|_V^{\Ga}]du\leq  \frac{\sup_{u\in[0,t]}|C(u)|+|x|_H^2}{\Gth}\cdot e^{c(t-s)}+\frac{\sup_{u\in[0,t]}|C(u)|}{\Gth}\cdot (t-s).
\end{align*}
Thus, by the Burkholder-Davis-Gundy inequality, there exists $C'>0$ such that
\begin{align*}
&\E[\sup_{s\leq u\leq t}|X(s,x;u)|_H^2]\\
&\leq C' \E\bigg[\int_{s}^{t}\bigg\{ \Vert B(t,X(s,x;u))\Vert_{L_2(H)}^2+\int_{\{|z|<1\}} |H(t,X(s,x;u),z)|_H^2 \nu(dz)\bigg\}du \bigg] \\
&\leq C' \E\bigg[\int_{s}^{t}\left\{\sup_{u\in[0,t]}|C(u)|+\Gg|X(s,x;u)|_V^{\Ga}+c|X(s,x;u)|_H^2\right\} du \bigg] \\
&\leq \frac{C'\left(\sup_{u\in[0,t]}|C(u)|+|x|_H^2\right)}{\Gth}\cdot\left[(\Gg+\Gth)(t-s)+(\Gg+c\Gth)e^{c(t-s)}\right],
\end{align*}
which implies that
\begin{eqnarray*}
\prob (\tau^x_R< t) &\leq& \prob (\sup_{u\in[s,t]}|X(s,x;u)|_H\geq R)\\
&\leq& \frac{\E[\sup_{u\in[s,t]}|X(s,x;u)|^2_H]}{R^2}\\
 &\leq& \frac{C'\left(\sup_{u\in[0,t]}|C(u)|+|x|_H^2\right)}{\Gth R^2}\cdot\left[(\Gg+\Gth)(t-s)+(\Gg+c\Gth)e^{c(t-s)}\right].
\end{eqnarray*}
Hence for $l\in(0,R-|x|_H)$, we have that
\begin{align}\label{AS2}
&\sup_{\lb w\in H:|w-x|_H\leq l \rb} \prob (\tau^w_R< t)\nonumber\\
&\ \ \leq \frac{C'\left(\sup_{u\in[0,t]}|C(u)|+[|x|_H+l]^2\right)}{\Gth R^2}\cdot\left[(\Gg+\Gth)(t-s)+(\Gg+c\Gth)e^{c(t-s)}\right].
\end{align}
Therefore, $P_{s,t}f$ is continuous at $x$ by (\ref{AS1}) and (\ref{AS2}). Since $x\in H$ is arbitrary, the proof for the case that $J\equiv 0$ is complete.

\vskip 0.2cm

We now consider the case that $J\not\equiv 0$. Let $\{Z(t)\}_{t\ge
0}$ be the unique solution of the SPDE (\ref{maineqnR}) with $J\equiv 0$.

\vskip 0.1cm

Denote by $P^Z(s,x;t,{\cal X})$ the transition semigroup of
$\{Z(t)\}_{t\ge
0}$,  where $x\in H, {\cal X}\in \sB(H)$ and $0\leq s< t<\infty$. Define $\zeta_1:= \inf\{u>s:N([s, u],\{|z|\ge1\})=1\}$, which is the
first jump time of $u\mapsto N([s,u],\{|z|\ge1\})$ after time $s$. Then, by conditioning on $\zeta_1$, we get
{\small\begin{align*}
P^X&(s,x;t,{\cal X})\\
&=e^{-\nu(\{|z|\ge1\})(t-s)}P^Z(s,x;t,{\cal X})\nonumber\\
&\ \ \ +
\int_{s}^t\int_{\{|x_2|\ge 1\}}
e^{-\nu(\{|z|\ge1\})(t_1-s)}P^X(t_1,x_1+J(t_1,x_1,x_2);t,{\cal X})P^Z(s,x;t_1,dx_1)\nu(dx_2)dt_1.
\end{align*}}

\vskip 0.1cm

Repeating this procedure, we get
\begin{eqnarray}\label{kl}
  P^X(s,x;t,{\cal X})=e^{-\nu(\{|z|\ge1\})(t-s)}\left[P^Z(s,x;t,{\cal X}) + \sum_{k=1}^{\infty}\Psi_k\right],
\end{eqnarray}
where
{\small$$
  \begin{array}{ll}
    \Psi_k= \displaystyle\mathop{\int\cdots\int}\limits_{s<t_1<\cdots<t_k<t}
   \int_{\{|x_2|\ge 1\}}\cdots\int_{\{|x_{2k}|\ge 1\}}
   {P}^{Z}(s,x;t_1,dx_1){P}^{Z}(t_1,x_1+J(t_1,x_1,x_2);t_2,dx_3)\\
   ~~~~~~~~~\times\cdots{P}^{Z}(t_k,x_{2k-1}+J(t_k,x_{2k-1},x_{2k});t,{\cal X})\nu(dx_2)\nu(dx_4)\cdots\nu(dx_{2k})dt_1dt_2\cdots dt_k.
\end{array}$$}  

Since we have shown that the transition
semigroup of $\{Z(t)\}_{t\ge
0}$ is strong Feller, $P^Z(s,x;t,{\cal X})$ and $\Psi_k$, $k\in \mathbb{N}$, are all continuous with
respect to $x$. Then, by \eqref{kl}, we conclude that $P^X(s,x;t,{\cal X})$ is lower semi-continuous with
respect to $x$. Therefore, the transition semigroup
$\{P_{s,t}\}$  is strong Feller by
\cite[Proposition 6.1.1]{MT}.
\end{proof}

\subsubsection{Proof of Theorem \ref{3-1}}
\begin{proof}

Denote the transition probability function of $(X(t),\GL(t))$ by $\lb P(s,(x,i),t,B\times \lb j\rb): 0\leq s<t, (x,i)\in H\times \cS, B\in \sB(H), j\in \cS\rb$. For $(x,i)\in H\times \cS$, let $X^{(i)}(t)$ be defined by Equation \eqref{RSThm1Eq1}. We also define  $\widetilde{X^{(i)}}(t)$ to be the killing process with generator $\sL+q_{ii}$. Then, for  $f \in B_b(H)$,
$$\E[f(\widetilde{X^{(i)}}(t))]=\E \bigg[ f(X^{(i)}(t))\exp\left\{ \int_0^{t} q_{ii}(X^{(i)}(u))du\right\}\bigg]. $$

\vskip 0.1cm

Let $\widetilde{P^{(i)}}(s,x;\cdot)$ be the transition probability function of $\widetilde{X^{(i)}}(t)$. Then, for $0\leq s<t$, $B\in \sB(H)$ and $j\in\cS$, we have
\begin{eqnarray*}
&&P(s,(x,i);t, B\times\lb j\rb)\\
& =& \delta_{ij} \widetilde{P^{(i)}}(s,x;t,B)\\
&& +\int_{s}^{t}\int_{H} P(t', (x',j'), t, B\times \lb j\rb)\bigg(\sum_{j'\in\cS\setminus\lb i\rb} q_{ij'}(x')\bigg)\widetilde{P^{(i)}}(s,x;t',dx')dt'.
\end{eqnarray*}

\noindent Repeating this procedure, we get
$$P(s,(x,i);t, B\times\lb j\rb)= \delta_{ij} \widetilde{P^{(i)}}(s,x;t,B)+\sum_{k=1}^{n}\Psi_k + U_n,$$
where
\begin{align*}
\Psi_k =& \int\cdots\int_{s<t_1<\cdots<t_k<t} \sum_{j_0,\dots,j_k}\int_{H^k}  q_{j_{k-1},j_{k}}(x_k)\widetilde{P^{(j_k)}}(t_k,x_k;t,B) \\
&\quad \times q_{j_{k-2},j_{k-1}}(x_{k-1})\widetilde{P^{(j_{k-1})}}(t_{k-1},x_{k-1};t_{k},dx_{k})\cdots   q_{i,j_{1}}(x_{1})\widetilde{P^{(j_1)}}(t_{1},x_{1};t_{2},dx_{2}) \\
&\quad \times \widetilde{P^{(i)}}(s,x;t_1,dx_1)dt_1\cdots dt_k
\end{align*}
and the sum is over
$$j_0=i,\ j_{\ell}\in \cS\setminus\lb j_{\ell-1}\rb \text{ for }\ell\in \lb 1,\dots,k-1\rb, \ j_k=j; $$
\begin{align*}
U_n =& \int\cdots\int_{s<t_1<\cdots<t_{n+1}<t} \sum_{j_0,\dots,j_{n+1}}\int_{H^{n+1}}  q_{j_{n},j_{n+1}}(x_{n+1})P(t_{n+1},x_{n+1};t,B\times \lb j\rb) \\
&\quad \times  q_{j_{n-1},j_{n}}(x_n)\widetilde{P^{(j_n)}}(t_n,x_n;t,B) \cdots   q_{i,j_{1}}(x_{1})\widetilde{P^{(j_1)}}(t_{1},x_{1};t_{2},dx_{2}) \\
&\quad \times \widetilde{P^{(i)}}(s,x;t_1,dx_1)dt_1\cdots dt_{n+1}
\end{align*}
and the sum is over
$$j_0=i,\ j_{\ell}\in \cS\setminus\lb j_{\ell-1}\rb \text{ for }\ell\in \lb 1,\dots,n+1\rb. $$

By condition \ref{Q0}, we find that $U_n\leq \frac{[(t-s)L]^{n+1}}{(n+1)!}$.
Letting $n\ra \infty$, we get
$$P(s,(x,i);t, B\times\lb j\rb)= \delta_{ij} \widetilde{P^{(i)}}(s,x;t,B)+\sum_{k=1}^{\infty}\Psi_k.$$
By Theorem \ref{3-2}, we know that the transition semigroup of $X^{(i)}(t)$ is strong Feller. Then, following the argument of \cite[Lemma 4.5]{XYZ}, we can show that the semigroup of $\widetilde{X^{(i)}}(t)$ is also strong Feller. Thus, we conclude that $\widetilde{P^{(i)}}(s,x;t,B)$ and $\Psi_k$ for $k\in\N$ are all continuous with  respect to $x$. Using the fact that $\cS$ is equipped with a discrete metric, we conclude that $P(s,(x,i);t, B\times\lb j\rb)$ is lower semi-continuous with respect to $(x,i)$. Therefore, the transition semigroup $\lb P_{s,t}\rb$ of $(X(t),\GL(t))$ is strong Feller by \cite[Proposition 6.1.1]{MT}.
\end{proof}

\subsection{Irreducibility}

For $\varpi>0$, define
$$
{D(A, \varpi)}:=\left\{ (v,i)\in V\times \cS: A(t,v,i)\in H\ {\rm and}\ \int_0^t|A(s,v,i)|_H^\varpi ds<\infty,\ \forall t \in [0,\infty)\right\}.
$$

\noindent We impose the following  assumptions.

\begin{enumerate}
\labitem{(D)}{D} There exists $\varpi >2$ such that $\overline{D(A,\varpi)}=H\times \cS.$

\labitem{(Q1)}{Q1} For any distinct $i,j
\in\mathcal{S}$, there exist an open set {$U\subset H$} and $j_1,\dots, j_r \in \mathcal{S}$ with
$j_p \neq j_{p+1}$, $j_1=i$ and $j_r=j$ such that $q_{j_{p}j_{p+1}}(x)>0$ for $p = 1,\ldots, r-1$ and $x\in U$.
\end{enumerate}

\vskip 0.1cm

\begin{assumption} \label{A4}  Assumption \ref{A1} holds with  $C\in L^{\infty}_{loc}([0,\infty);(0,\infty))$, $\gamma=0$ and the exponent  $\alpha$ in condition \ref{growthRho} replaced by some $\Ga'\in (1,\Ga)$.

\end{assumption}

\vskip 0.1cm

\vskip 0.2cm

Now we give the main result of this subsection.
\vskip 0.4cm

\begin{theorem}\label{4-1}
Suppose that Assumption \ref{A4}  and conditions \ref {D}, \ref{Q1} hold. Then, the transition semigroup $\lb P_{s,t}\rb$ of $(X(t),\GL(t))$ is irreducible.
\end{theorem}

\vskip 0.1cm

\begin{proof}
First, we consider the case that $i\in\cS$ is fixed. To simplify notation, we drop the dependence on $i$. Define the
first jump time of $\{X(t)\}_{t\ge 0}$ by
$$\zeta_1:= \inf\{t>0:N([0, t],\{|z|\ge1\})=1\},$$
which is exponentially distributed with rate $\nu(\{|z|\ge1\})$. Let $\{Z(t)\}_{t\ge
0}$ be the unique solution of the SPDE (\ref{maineqnR}) with $J\equiv 0$. We have that
$Z(t)=X(t)$ for $t<\zeta_1$. Hence, to obtain the irreducibility of $\{P_{s,t}\}$, we may assume without loss of generality that $J\equiv 0$.

\vskip 0.1cm

Denote by $\{X^x(t)\}_{t\ge 0}$ the solution of Equation \eqref{maineqnR} with $X(0)=x$, $x\in H$. Let $T,M,R>0$, $t_1\in (0,T)$, $y(0)\in D(A_H)$ and $\{y(t)\}_{t\in[0,T]}$ be the solution of the following equation:
\begin{align}\label{auxeqn1}
dy(t) &= A(t,y(t))dt-\frac{M}{T-t_1} (y(t)-y(0))dt,\ \ t> t_1, \nonumber\\
y(t)&= X^x(t_1)1_{\lb |X^x(t_1)|_H\leq R\rb},\ \ t\le t_1.
\end{align}

\vskip 0.1cm

Similar to \cite[Lemma 2.1]{Z}, we can prove the following result.

\vskip 0.1cm

\begin{lemma}\label{auxlemma1}
Let $m:=2M-(K+\rho(y(0))+1)T> 0$. Then,
\begin{align}\label{auxest1}
|y(t)-y(0)|_H^2\leq e^{-\frac{m (t-t_1)}{T-t_1}}(R+|y(0)|_H)^2+\int_{t_1}^{t}e^{-\frac{m (t-s)}{T-t_1}}  |A(s,y(0))|_H^2ds,\ \ t\in [t_1,T],
\end{align}
\begin{align*}
\int_{t_1}^{T}|y(t)-y(0)|_H^2dt&\leq \frac{T-t_1}{m}\left[ (R+|y(0)|_H)^2+\int_{t_1}^{T}  |A(s,y(0))|_H^2ds\right],
\end{align*}
and there exists $\vartheta>0$, which is independent of $t_1$,  such that
\begin{align*}
\int_{t_1}^{T}\rho(y(s))ds\leq  \vartheta(T-t_1)^{\frac{\Ga-\Ga'}{\Ga}}.
\end{align*}
\end{lemma}

\vskip 0.1cm

Note that in this paper $A$ is time-dependent. Hence we need replace $(T-t_1)^2|A(y_0)|^2$ in \cite[Lemma 2.1]{Z} by $\int_{[-t_1,T]}w(s)|A(s,y(0))|^2ds$ with suitable $w(s)$.

\vskip 0.1cm

Let $\Ge>0$ and $n\in \N$. We consider the following equation:
\begin{align} \label{auxeqn2}
d\widetilde{X}^n(t)&=A(t,\widetilde{X}^n(t))dt+B(t,\widetilde{X}^n(t))dW(t)+\int_{\{|z|<1\}} H(t,\widetilde{X}^n(t),z)\widetilde{N}(dt,dz) \nonumber\\
&\quad -\frac{M}{T-t_1}(\Ge B_n^{-1}+I)^{-1}(y(t)-y(0))\chi_{[t_1,T]}(t)dt, \nonumber\\
\widetilde{X}^n(0)&=x,
\end{align}
where $\{y(t)\}_{t\in[0,T]}$ is the solution of \eqref{auxeqn1}. Note that $\widetilde{X}^n(t_1)=X(t_1)$.
By Theorem \ref{2-2} and \eqref{auxest1}, we know that \eqref{auxeqn2} has a unique solution $\{\widetilde{X}^n(t)\}_{t\in[0,T]}$.

\vskip 0.1cm

Following the argument of \cite{Z} and carefully handling the dependence of constants, we obtain the following estimation for the solution $\widetilde{X}^n$.

\vskip 0.1cm

\begin{lemma}\label{auxestimatelemma}
There exists $\vartheta>0$, which is independent of $\Ge,n,T,M,R,t_1$,  such that
\begin{align*}
&\sup_{n\in\N}\bigg\{\E\bigg[\sup_{s\in[0,T]}|\widetilde{X}^n(s)|_H^2\bigg]\bigg\}\nonumber \\
\leq& e^{\vartheta T}\bigg(\vartheta(T+|x|^2_H)+\frac{M^2}{m(T-t_1)} \bigg[(R+|y(0)|_H)^2+\int_{t_1}^{T}|A(s,y(0))|_H^2ds\bigg]\bigg).
\end{align*}
\end{lemma}

\vskip 0.1cm

By virtue of Lemma \ref{auxestimatelemma}, following the argument of \cite[Lemma 2.3]{Z}, we can prove the following lemma.
\begin{lemma}\label{lemmaIRR}
Let $T>0$, $y(0)\in D(A,\varpi)$ and $\eta,\Gd\in (0,1)$. Then, there exist $M,R>0$ and $t_1\in (0,T)$ such that for any $n\in \N$ we can find an $\Ge\in (0,1)$ satisfying
$$\prob(|\widetilde{X}^n(T)-y(0)|_H>\Gd)<\eta. $$
\end{lemma}

By Lemma \ref{lemmaIRR} and following the argument of \cite[Theorem 1.1]{Z}, we complete the proof of Theorem \ref{4-1} for the case that $i\in\cS$ is fixed. Finally, similar to \cite[Theorem 3]{GS2}, we can complete the proof for the general case.
\end{proof}

\section{Existence and Uniqueness of Periodic Measures}

In this section, we discuss periodic measures for the hybrid system (\ref{maineqnR}) and (\ref{222}) .

\begin{definition}
Let  $E$ be a Polish space with Borel $\sigma$-algebra $\sB (E)$,  $\lb Y(t)\rb_{t\geq 0}$ an $E$-valued Markov process, and $\ell>0$. A  probability measure $\mu_0$ on $\sB (E)$
 is said to be an $\ell$-periodic measure for $\lb Y(t)\rb_{t\geq 0}$   if the following condition holds:

\begin{itemize}
\item $Y(0)$ has distribution $\mu_0$ implies that the joint distribution of $Y(t_1+k\ell),\dots, Y(t_n+k\ell)$ is independent of $k$ for all $0\leq t_1<\cdots<t_n$ and $n\in\N$.
\end{itemize}

\end{definition}

 For $g\in
C^{1,2}([0,\infty)\times H\times\mathcal{S};\mathbb{R})$,
define
\begin{eqnarray*}
{\sA}g(t,x,i):={\sL}_ig(t,x,i)+Q(x)g(t,x,\cdot)(i)
\end{eqnarray*}
with
\begin{eqnarray*}
&&{\sL}_ig(\cdot,\cdot,i)(t,x)\\&:=&g_t(t,x,i)+\langle A(t,x,i),g_x(t,x,i) \rangle+
\frac{1}{2}\mathrm{ trace}(B^T(t,x,i)g_{xx}(t,x,i)B(t,x,i))\nonumber\\
&&+\int_{\{|z|<1\}}[g(t,x+H(t,x,i,z),i)-g(t,x,i)-
\langle g_x(t,x,i), H(t,x,i,z)\rangle]\nu(\mathrm{d}z)\nonumber\\
&&+\int_{\{|z|\ge 1\}}[g(t,x+J(t,x,i,z),i)-g(t,x,i)]\nu(\mathrm{d}z),
\end{eqnarray*}
and
\begin{eqnarray*}
Q(x)g(t,x,\cdot)(i):=\sum_{j\in\mathcal{S}}[g(t,x,j)-g(t,x,i)]q_{ij}(x).
\end{eqnarray*}

\vskip 0.1cm

We make the following assumption  on the matrix $Q=(q_{ij}(x))$.

\vskip 0.1cm

\begin{enumerate}
\labitem{(Q2)}{Q2} There exists a positive increasing function $f$ on ${\cal S}$ satisfying
 $$
 \lim_{j\rightarrow\infty }f(j)=\infty,\ \ \
 \sup\limits_{x\in H,\,i\in\mathcal{S}}\sum\limits_{j\neq i} [f(j)-f(i)]q_{ij}(x) <\infty,\ \ \ \lim_{i\rightarrow\infty }\sup\limits_{x\in H}\sum\limits_{j\neq i} [f(j)-f(i)]q_{ij}(x)=-\infty.$$
\end{enumerate}

\vskip 0.1cm

Now we can state the main theorem of this paper.

\begin{theorem} \label{existencePSRS}
Let $\ell>0$. Suppose that functions $A, B,H,J$ are all $\ell$-periodic with respect to $t$, the embedding of $V$ into $H$ is compact, Assumption \ref{A4} holds with $C\in L^{\infty}([0,\infty);(0,\infty))$, Assumption \ref{A2} and conditions \ref{D}, \ref{Q0}, \ref{Q1}, \ref{Q2} hold, and
\begin{eqnarray}\label{Cond6}
&&\lim_{n\ra \infty} \left.\sup_{|v|_V>n, t\ge0,i\in\cS} \left\{-\theta |v|_V^{\Ga}+c|v|_H^2\right.\right.\nonumber\\
&&\ \ \ \ \ \ \ \ \ \ \ \ \ \ \ +\left.\int_{\{|z|\ge1\}}\left [|J(t,v,i,z)|_H^2+2\lan v, J(t,v,i,z) \ran_H\right]\nu(dz)\right\}\nonumber\\
&=&-\infty.
\end{eqnarray}

\vskip 0.1cm

\noindent Then,

\vskip 0.1cm

 (i) Equations \eqref{maineqnR} and \eqref{222} have a unique solution $\{(X(t),\GL(t))\}_{t\geq 0}$;

\vskip 0.1cm

 (ii) The transition semigroup $\{P_{s,t}\}$ of $\{(X(t),\GL(t))\}_{t\geq 0}$ is strong Feller and irreducible;

\vskip 0.1cm

(iii) The hybrid system $\{(X(t),\GL(t))\}_{t\geq 0}$ has a unique $\ell$-periodic measure $\mu_0$;

\vskip 0.1cm

 (iv) Let $\mu_s(A)=\mathbb{P}_{\mu_0}((X(s),\GL(s))\in A)$ for $A\in
\mathcal{B}(H\times \cS)$ and $s\ge 0$. Then,
for any
$s\ge 0$ and $\varphi\in L^2(H\times \cS;\mu_s)$,
\begin{equation}\label{May14v}
\lim_{n\rightarrow\infty}
\frac{1}{n}\sum_{i=1}^{n}P_{s,s+i\ell}\varphi=\int_{H\times \cS}\varphi
d\mu_s\ \
{\rm in}\ \ L^2(H\times \cS;\mu_s).
\end{equation}

\end{theorem}


\vskip 0.2cm

\begin{proof}

Claims  \textit{(i)} and  \textit{(ii)} follow from Theorems \ref{2-1},  \ref{3-1} and \ref{4-1}.


\vskip 0.1cm

Let $(X(t),\GL(t))$ be the unique solution to the hybrid system  \eqref{maineqnR} and \eqref{222} with initial value $(x,i)\in H\times \cS$. Define
$\cV(t,x,i)=|x|_H^2+f(i)$, where $f$ is given by condition  \ref{Q2}.
Then,
$$\sA {\cal V} (t,x,i)=\sL_i |x|_H^2 +\sum_{j\in\cS} [f(j)-f(i)]q_{ij}(x),$$
where
\begin{align*}
\sL_i |x|_H^2 &=  2\lan A(t,x,i), x\ran+\frac{1}{2}\|B(t,x,i)\|_{L_2(H)}^2+\int_{\{|z|<1\}} |H(t,x,i,z)|_H^2 \nu(dz) \\
&\ \ \ +  \int_{\{|z|\ge1\}}\left [|J(t,x,i,z)|_H^2+2\lan x, J(t,x,i,z) \ran_H\right]\nu(dz).
\end{align*}

\noindent Then, by conditions \ref{coercivityr},  \eqref{Cond6} and \ref{Q2}, we obtain that
\begin{align}
&\lim_{|y|_H+i\ra \infty} \inf_{t\ge0}{\cal V} (t,y,i)=\infty, \nonumber\\
&\label{V4C2}\lim_{n\ra \infty} \sup_{|y|_V+i>n, i\in\cS, t\ge0} \sA {\cal V} (t,y,i)=-\infty,\\
&\sup_{y\in V, i\in\cS, t\ge0} \sA {\cal V} (t,y,i) <+\infty.\nonumber
\end{align}

For $n\in\N$, define the stopping time $T_n$ by
$$T_n:=\inf\lb t\geq 0: |X(t)|_{V}\vee \GL(t) \geq n\rb. $$
For $t\ge0$, by It\^o's formula (see Gy\"{o}ngy and  Krylov \cite[Theorem 2]{GK}), we get
\begin{align} \label{25A2}
&\E[{\cal V}(t\wedge T_n, X(t\wedge T_n), \GL(t\wedge T_n))] \nonumber\\
&= \E[{\cal V}(0,X(0),\GL(0))]+\E\bigg[\int_0^{t\wedge T_n} \sA{\cal V}(u,X(u),\GL(u))du\bigg].
\end{align}
Define
$$A_n:= -\sup_{|y|_V+k>n, t\ge0} \sA{\cal V}(t,y,k). $$
By \eqref{V4C2}, we get $\ds \lim_{n\ra \infty} A_n=\infty$.

\vskip 0.1cm

We have
$$\sA{\cal V}(u,X(u),\GL(u)) \leq -1_{\lb |X(u)|_V+k\geq n\rb} A_n+\sup_{|y|_V+k<n, u\ge0} \sA{\cal V}(u,y,k). $$
Then, there exist positive constants $\epsilon_1$ and $\epsilon_2$ such that for sufficiently large $n$,
\begin{align}\label{c1c2n0}
\E\bigg[\int_0^{t\wedge T_n} 1_{\lb |X(u)|_V+k \geq n\rb} du \bigg] &\leq \frac{\epsilon_1 t +\epsilon_2}{A_n}.
\end{align}
Letting $n\ra \infty$ in \eqref{c1c2n0}, we get
\begin{align} \label{25E'}
\lim_{n\ra\infty} \limsup_{T\ra\infty} \frac{1}{T} \int_0^T P(0,(x,i);u, B_n^c)du =0,
\end{align}
where $B_n^c= \lb (y,k)\in H \times \cS: |y|_V+k \geq n\rb$.

\vskip 0.1cm

By (\ref{25A2}), we find that there exists $\lambda>0$ such that
$$\mathbb{E}[{\cal V}(s,X(s),\GL(s))]\leq \lambda s+{\cal V}(0,x,i),\ \ s>0,$$
which together  with Chebyshev's inequality implies that
$$
\mathbb{P}(0,(x,i);s,B_n^c)\leq \frac{\lambda s+{\cal V}(0,x,i)}{\inf_{|y|_V> n,\, i\in \cS, t>0}{\cal V}(t,y,i)}.
$$
Hence, there exists a sequence of positive integers $\gamma_n\uparrow\infty$ such that
\begin{equation}\label{25F'}
\lim_{n\rightarrow\infty}\left\{\sup_{(x,i)\in
B_{H\times \cS}(\gamma_n),\,s\in(0,\ell)}\mathbb{P}(0,(x,i);s, B_n^c)\right\}=0,
\end{equation}
where $B_{H\times \cS}(\iota):=\{(y,k)\in H\times \cS:|y|_H+k<\iota\}$ for $\iota>0$.


\vskip 0.1cm

By the assumption that functions $A, B,H,J$ are all $\ell$-periodic with respect to $t$, we find that the transition semigroup $\{P_{s,t}\}$ is $\ell$-periodic, i.e.,
$$P(s,(x,i);t,A)= P(s+\ell,(x,i);t+\ell,A),\ \ \ \ \forall 0\le s<t, x\in H, \ i\in \cS, \  A\in \sB(H\times \cS).
 $$
Since the embedding of $V$ into $H$ is compact, combining the periodicity and the strong Feller property of $\{P_{s,t}\}$ with
(\ref{25E'}), (\ref{25F'}) and following the argument of \cite[Theorem 3.2 and Remark
3.1]{K}, we conclude that $\{(X(t),\GL(t))\}_{t\geq 0}$ has an $\ell$-periodic measure $\mu_0$.

\vskip 0.1cm

By using the same argument, we can show that \cite[Lemma 3.12 and Theorem
3.13]{GS} hold with the state  space $\mathbb{R}^m$ replaced by
$H\times\cal{S}$. Then, by Theorems \ref{3-1} and \ref{4-1}, we conclude that
there exists a unique family of probability measure  $\lb\eta_s\rb$ on $\sB(H\times \cS)$
that is $\ell$-periodic with respect to $\lb P_{s,t} \rb$, that is,
  $$\eta_s(A)=\int_H P(s,x; s+\ell, A)\eta_s(dx),\ \ \ \ \forall s\ge0, A\in \sB(H\times \cS).$$
Hence, we obtain the uniqueness of periodic measures, namely, $\mu_0:= \eta_0$. Finally, following the same argument of the proof of \cite[Theorem 3.13]{GS}, we obtain (\ref{May14v}).
\end{proof}

\section{An Example}

\begin{example*}[Stochastic porous media equations]\label{26H}
Let $d\geq 1$ and $\cO\subset \R^d$ be a bounded open set with smooth boundary. Set $\mathfrak{L}=-(-\GD)^{\Gg}$ for $\Gg>0$. Let $r>1$. Denote by $d\mu$ the normalized Lebesgue measure on $\cO$ and define
$$
V=L^{r+1}\left(\cO; d\mu\right),\ \ \ \ \ H=H^{-\Gg}\left(\cO;d\mu\right),
$$
where
$H^{-\Gg}(\cO; d\mu) $ is the completion of $L^2(\cO;d\mu)$ with respect to the norm
$$|f|_{H^{-\Gg}}= \left(\int_{\cO}|(-\GD)^{-\frac{\Gg}{2}}f|^2d\mu\right)^{\frac{1}{2}},\ \ f\in L^2\left(\cO; d\mu\right).$$
Note that the embedding of $V$ into $H$ is compact. 

\vskip 0.1cm

Let $\lb q_{ij}\rb$ be  measurable functions defined on $H$ such that one of the following conditions is satisfied:

\vskip 0.1cm

(a) There exist $m\in \N$ and $M>0$ such that for all $x\in H$,
\begin{align*}
q_{ij}(x) &=0  \quad \text{if }|i-j|>m, \\
q_{ij}(x) &\in (0,M] \quad \text{if }0<|i-j|\leq m,
\end{align*}
and
\begin{align*}
\inf_{x\in H, i>m, j\in[i-m,i)} q_{ij}(x) &> \sup_{x\in H, i>m, j\in(i,i+m]} q_{ij}(x).
\end{align*}

\vskip 0.1cm

(b) $$
0<\inf_{x\in H,\,j\neq i}\{j^{1+\delta}q_{ij}(x)\}< \sup_{x\in H,\,j\neq  i}\{j^{1+\delta}q_{ij}(x)\}<\infty \ \ {\rm for\ some}\ \delta>0.
$$

\noindent Then, conditions \ref{Q0},  \ref{Q1} and \ref{Q2} (setting $f(i)=i^2$ under condition (a), and $f(i)=i^{\delta/2}$ under condition  (b)) hold (cf. \cite[Examples 1 and 2]{GS2}).

\vskip 0.1cm

Let $\cS=\N$, $\ell>0$, and $g$ be a measurable function on  $[0,\infty)\times \cS$ such that $g(\cdot, i)\in L^{\infty}[0,\infty)$ is $\ell$-periodic for each $i\in \cS$ and $\sup_{t\geq 0, i\in\cS}|g(t,i)|<+\infty$. For $t\in [0,\infty)$ and $x\in \R$, define
$$
\Psi(x)=|x|^{r-1}x,\ \ \ \ \Phi(t,x,i)= g(t,i)x,
$$
and
$$
A(t,x,i)=\kappa(t,i)\mathfrak{L}(\Psi(x))+\Phi(t,x,i),
$$
where $\kappa(\cdot, i)$ is $\ell$-periodic for each $i\in \cS$ and $k_1\leq\kappa(t,i)\leq k_2$ for all $t\in [0,\infty)$, $i\in{\cal S}$  and some positive constants $k_1<k_2$.
One can show that condition \ref{D} holds and there exists $C_r>0$ such that for all $t\in [0,\infty)$, $x_1,x_2\in V$  and $i\in\cS$,
\begin{align} \label{5A}
\lan A(t,x_1,i)-A(t,x_2,i),x_1-x_2\ran \leq -\kappa(t,i)C_{r}|x_1-x_2|_V^{r+1}+g(t,i)|x_1-x_2|_H^2.
\end{align}

\vskip 0.1cm

Let $\lambda_1\le\lambda_2\le\cdots\le\Gl_j\leq\cdots$ be the eigenvalues of $-\GD$ and $\{e_j \}$ the corresponding unit eigenvectors. Suppose $\frac{1}{2}<s\leq \frac{\Gg}{d}$. Define $B_0\in L_2(H)$ by
$$
B_0e_j=j^{-s}e_j,\ \ j\in \mathbb{N}.
$$
Let $b>0$ and $\lb b_j\rb_{j\in \mathbb{N}}$ be measurable functions defined  on  $[0,\infty)\times H\times \cS$ satisfying the following conditions:
\begin{eqnarray*}
&&b_j(t,x,i)=b_j(t+\ell,x,i),\ \ t\in[0,\infty),x\in H,i\in\cS,j\in\mathbb{N},\\
&&|b_j(t,x,i)-b_j(t,y,i)|\leq b|x-y|_{H},\ \ t\in[0,\ell),x,y\in H,i\in\cS,j\in\mathbb{N},\\
&&\sup_{t\in[0,\ell),x\in H,i\in \cS,j\in\mathbb{N}} |b_j(t,x,i)|\le b,\\
&&\inf_{t\in[0,\ell),|x|_H\leq n,i\in \cS,j\in\mathbb{N}} b_j(t,x,i)>0,\ \ \ \ \forall n\in \N.
\end{eqnarray*}
Define
$$
B(t,x,i)e_j=b_j(t,x,i)j^{-s}e_j,\ \ t\in[0,\infty),x\in H,i\in\cS,j\in\mathbb{N}.
$$
As an explicit example, similar to \cite[Example 4.6]{Z}, we may let
$$
b_j(t,x,i):= \frac{b'(t,i)}{1+j^{-\frac{2\gamma}{d}}|\lan x,e_j\ran_{L^2(\cO)}|},
$$
where $b'$ is a measurable function on  $[0,\infty)\times \cS$ such that $b'(\cdot, i)$ is $\ell$-periodic for each $i\in \cS$ and $0<\inf_{t\in[0,\ell),i\in \cS}b'(t,i)\le \sup_{t\in[0,\ell),i\in \cS}b'(t,i)<\infty$.

\vskip 0.1cm

Let $\{W(t)\}_{t\ge 0}$ be an $H$-valued cylindrical Wiener process on a complete filtered probability space $(\GO, \sF, \{\sF_t\}_{t\geq 0}, \prob)$, $Z$ a real Banach space and $N$ a Poisson random measure on $(Z,{\sB}(Z))$
with intensity measure $\nu$. Assume that $W$ and $N$ are independent.
Suppose that $c>0$, $K>0$,  $\rho\in L^{\infty}_{loc}(V;[0,\infty))$ and $H,J: [0,\infty)\times V\times \cS\times Z \ra H$ are measurable functions satisfying the following conditions:
\begin{eqnarray*}
&&H(t,x,i,z)=H(t+\ell,x,i,z),\ \ \ \ J(t,x,i,z)=J(t+\ell,x,i,z),\ \ t\geq 0,x\in V,i\in\cS,z\in Z,\\
&&\int_{\{|z|<1\}} |H(t,x_1,i,z)-H(t,x_2,i,z)|_H^2 \nu(dz) \leq (K+\rho(x_2))|x_1-x_2|_H^2,\ \ t\ge0,x_1,x_2\in V, i\in \cS,\\
&&\int_{\{|z|<1\}} |H(t,x,i,z)|_H^2 \nu(dz) \leq c(1+|x|_H^2),\ \ \ \ \int_{\{|z|\ge1\}}|J(t,x,i,z)|^2_H \nu(dz)\le c,\ \ t\geq 0,x\in V,i\in\cS.
\end{eqnarray*}

\vskip 0.1cm

Define
$$
\ds B_n:= \bigg(\inf_{t\in[0,\ell),|x|_H\leq n,i\in \cS,j\in\mathbb{N}} b_j(t,x,i)\bigg) B_0,\ \ n\in \mathbb{N}.
$$
We now show that all conditions of Theorem \ref{existencePSRS} are satisfied. 
It is easy to see that conditions \ref{hemicontr} and \ref{nondegenr} hold. Note that
\begin{align*}
&2\lan A(t,x,i),x\ran+ \Vert B(t,x,i)\Vert_{L_2(H)}^2 +\int_{\{|z|<1\}} |H(t,x,i,z)|_H^2 \nu(dz) \\
&\leq  - 2\kappa(t,i)|x|_V^{r+1}+[2g(t,i)+c]|x|_H^2 +\sum_{j=1}^{\infty}\frac{b^2}{j^{2s}}+  c.
\end{align*}
Hence, condition \ref{coercivityr} is satisfied by taking $C(t)=\sum_{j=1}^{\infty}\frac{b^2}{j^{2s}}+c$, $\Gth=2k_1$ and $\alpha=r+1$.
We have that
\begin{align*}
|A(t,v,i)|_{V^*} &\leq  k_2|-(-\GD)^{\Gg}(\Psi(v))|_{V^*}+ |g(t,i)|\cdot|v|_V \\
& \leq k_2(1+ c_{r} |v|_{V}^{r})+|g(t,i)| (c_{r}'+c_{r}''|v|_V^{r})\\
&\le \left(k_2+c_{r}'\sup_{t\geq 0, i\in\cS}|g(t,i)|\right)+\left(k_2c_{r}+c_{r}''\sup_{t\geq 0, i\in\cS}|g(t,i)|\right)|v|_V^{r},
\end{align*}
where $c_r, c_r', c_r''$ are the constants from Young's inequality for products. Then, we get \ref{growthAr} by modifying $C(t)$ to be the maximum of $\sum_{j=1}^{\infty}\frac{b^2}{j^{2s}}+c$ and
$$
2^{1/r} \left(k_2+c_{r}'\sup_{t\geq 0, i\in\cS}|g(t,i)|\right)^{\frac{r+1}{r}}.
$$
Additionally, set $\Gb =0$.

\vskip 0.1cm

Conditions \ref{growthBHr} and \ref{growthHr} follow from the definitions of  $B$ and $H$.  By \eqref{5A} and \cite[Theorem 2.4.1]{W}, we get
\begin{align*}
&\lan A(t,x_1,i)-A(t,x_2,i),x_1-x_2\ran\\
&\leq -\kappa(t,i)C_{r}|x_1-x_2|_V^{r+1}+g(t,i)|x_1-x_2|_H^2\\
&\leq -\kappa(t,i)C_r |B_n^{-1}(x_1-x_2)|_H^{\Gl}|x_1-x_2|_H^{r+1-\Gl}+g(t,i)|x_1-x_2|_H^2.
\end{align*}
Then  \ref{Localmonotone2r} and  hence \ref{localmonr} holds.
Finally,  by the assumption on $J$ and $\alpha=r+1>2$, we find that condition \eqref{Cond6} holds.  Thus, all conditions of Theorem \ref{existencePSRS} are fulfilled. Therefore, all assertions of Theorem \ref{existencePSRS} hold.
\end{example*}

\vskip 0.5cm

\begin{large} \noindent\textbf{Acknowledgements} \end{large} This work was supported by the Natural Sciences and Engineering Research Council of Canada (No. 4394-2018).


\vskip 0.5cm
\begin{thebibliography}{1234}

\bibitem{A} D. Applebaum. L\'evy Processes and Stochastic Calculus,
Second Edition. Cambridge University Press, (2009).

\bibitem{ATW}  M. Arnaudon, A. Thalmaier and F.Y. Wang. Harnack inequality and heat kernel estimates on
manifolds with curvature unbounded below. Bull. Sci. Math. 130, 223–233 (2006).


\bibitem{BLZ} Z. Brze\'zniak, W. Liu and J. Zhu. Strong solutions for SPDE with locally monotone coefficients driven by L\'evy noise. Nonl. Anal.: Real World Appl. 17, 283-310 (2014).



\bibitem {CHLY}  F. Chen, Y. Han, Y. Li and X. Yang. Periodic solutions of Fokker-Planck equations. J. Diff. Equ. 263, 285–298 (2017).


\bibitem{CL} M. Cheng and Z. Liu. Periodic, almost periodic and almost automorphic solutions for SPDEs with monotone coefficients. Disc. Contin. Dyn. Syst., Ser. B. 26, 6425-6462 (2021).


\bibitem{DD} G. Da Prato and A. Debussche. 2D stochastic Navier–Stokes equations with
a time-periodic forcing term. J. Dyn. Diff. Equ. 20, 301-335 (2008).

\bibitem{DT} G. Da Prato and C. Tudor. Periodic and almost periodic solutions for semilinear stochastic equations. Stoch. Anal. Appl. 13, 13–33 (1995).

\bibitem{DZ} G. Da Prato and J. Zabczyk.
Ergodicity for Infinite-dimensional Systems. Cambridge University Press (1996).

\bibitem{DZ3} G. Da Prato and J. Zabczyk.
Stochastic Equations in Infinite Dimensions, Second Edition. Cambridge University Press (2014).






\bibitem{FZ} C. Feng and H. Zhao. Random periodic solutions of SPDEs via integral
equations and Wiener–Sobolev compact embedding. J. Funct. Anal. 262, 4377–4422 (2012).

\bibitem{GM} L. Gawarecki and V. Mandrekar. Stochastic Differential Equations in Infinite Dimensions with Applications to Stochastic Partial Differential Equations. Springer, (2011).


\bibitem{GR} B. Gess and M. R\"{o}ckner. Stochastic variational inequalities and regularity for degenerate stochastic partial differential equations. Trans. Amer. Math. Soc. 369,  3017-3045 (2017).

\bibitem{GS} X.X. Guo and W. Sun. Periodic solutions of stochastic differential equations driven by L\'evy
noises. J. Nonl. Sci. 31:32  (2021).

\bibitem{GS2} X.X. Guo and W. Sun. Periodic solutions of hybrid jump diffusion processes. Front. Math. China 16, 705–725 (2021).

\bibitem{GK} I. Gy\"{o}ngy and N.V. Krylov. On stochastic equations with respect to semimartingales II, Ito formula in Banach spaces. Stoch. 6, 153-174 (1982).

\bibitem{HM} M. Hairer and J. C. Mattingly. Ergodicity of the 2D Navier-Stokes equations with degenerate stochastic
forcing. Ann. Math. 164, 993–1032 (2006).

\bibitem{a12} H. Hu and L. Xu. Existence and uniqueness theorems for periodic
Markov process and applications to stochastic functional differential
equations. J. Math. Anal. Appl. 466, 896-926 (2018).



\bibitem{JQSY} M. Ji, W. Qi, Z. Shen and Y. Yi. Existence of periodic probability solutions to Fokker-Planck equations with applications, J. Funct. Anal. 277,  108281 (2019).

\bibitem{K} R.Z. Khasminskii. Stochastic Stability of
Differential Equations, Second Edition.  Springer (2012).




\bibitem{LR} W. Liu and W. R\"{o}ckner. Stochastic Partial Differential Equations: an Introduction. Springer (2015).


\bibitem{MS} B. Maslowski and J. Seidler.  Invariant measures for nonlinear SPDE's: uniqueness and stability. Arch. Math. 34,  153-172 (1998).

\bibitem{M} M. M\'etivier. Semimartingales. A Course on Stochastic Processes. De Gruyter, (1982).

\bibitem{MT} S.P. Meyn and R.L. Tweedie. Markov Chains and Stochastic Stability. Springer, (1993).

\bibitem{Neuss} M. Neu\ss. Ergodicity for singular-degenerate stochastic porous media equations. J. Dyn. Diff. Equ. (2021). https://doi.org/10.1007/s10884-021-09961-9



\bibitem{PZ} S. Peszat and J. Zabczyk.
  Stochastic Partial Differential Equations with L\'{e}vy Noise, an Evolution Equation Approach.
 Cambridge University Press (2007).

\bibitem{RX} M. Romito and L. Xu. Ergodicity of the 3D stochastic Navier–Stokes equations driven by mildly degenerate noise. Stoch. Proc. Appl. 121, 673–700  (2011).

\bibitem{W} F.Y. Wang. Harnack Inequalities for Stochastic Partial Differential Equations. Springer (2013).

\bibitem{XYZ} F. Xi, G. Yin and C. Zhu. Regime-switching jump diﬀusions with non-Lipschitz coefficients and countably many switching states: existence and uniqueness, Feller, and strong Feller properties. Modeling, Stochastic Control, Optimization, and Applications. 571-599, (2019).

\bibitem{X} B. Xie. Uniqueness of invariant measures of infinite dimensional stochastic differential equations driven by Lévy noises. Potent. Anal. 36, 35-66 (2012).

\bibitem{a9} D. Xu D., Y. Huang and Z. Yang. Existence theorems
for periodic Markov process and stochastic functional differential
equations. Disc. Contin. Dyn. Syst. Ser. A. 24, 1005-1023 (2009).

\bibitem{ra} D. Xu, B. Li, S. Long and L. Teng. Moment
estimate and existence for solutions of stochastic functional differential
equations. Nonl. Anal. 108, 128-143 (2014).





\bibitem{YB} C. Yuan and J. Bao. On the exponential stability of switching-diffusion processes with jumps. Quart. Appl. Math. 71, 311–329 (2013).

\bibitem{Zr} R. Zhang. Existence and uniqueness of invariant measures of 3D stochastic MHD-$\Ga$ model driven by degenerate noise. Appl.  Anal. 101, 629–654 (2020).


\bibitem{Z} S.Q. Zhang. Irreducibility and strong Feller property for non-linear SPDEs. Stoch. 91, 352-382 (2019).

\bibitem{ZWL} X. Zhang, K. Wang and D. Li. Stochastic periodic solutions of stochastic differential equations driven by L\'{e}vy process. J. Math. Anal. Appl. 430, 231-242 (2015).

\end{thebibliography}
\end{document}